\documentclass[11pt]{amsart}

\usepackage{amsmath,amsthm, amscd, amssymb, amsfonts}
\usepackage[all]{xy}
\usepackage{enumitem}
\usepackage{array}  
\usepackage{multicol}
\usepackage{nicefrac}
\newcolumntype{L}{>{$}l<{$}} 
\newcolumntype{R}{>{$}r<{$}}
\newcolumntype{C}{>{$}c<{$}}
\usepackage{hyperref}
\usepackage{mathrsfs}

\usepackage{fontsmpl}
\usepackage{srcltx}
\usepackage[toc,page]{appendix}

\usepackage[T2A,T1]{fontenc}

\allowdisplaybreaks

\newcommand{\ba}{\mathbf{a} }

\newcommand{\cI}{{\mathcal I}}
\newcommand{\cJ}{{\mathcal J}}

\newcommand{\cZ}{\mathcal{Z}}

\newcommand{\toba}{\mathscr{B}}
\newcommand{\wtoba}{\widetilde{\toba}}
\newcommand{\htoba}{\widehat{\toba}}
\newcommand{\tobaq}{\toba_{\bq}}
\newcommand{\tobaqsr}{\breve{\toba}}

\newcommand{\Dchaintwo}[3]{\xymatrix@C-4pt{\overset{#1}{\underset{\ }{\circ}}\ar
@{-}[r]^{#2}
& \overset{#3}{\underset{\ }{\circ}}}}
\newcommand{\Dchainthree}[5]{\xymatrix@C-6pt{
\overset{#1}{\underset{\ }{\circ}}\ar  @{-}[r]^{#2}  & \overset{#3}{\underset{\
}{\circ}}\ar  @{-}[r]^{#4}
& \overset{#5}{\underset{\ }{\circ}} }}

\newcommand{\Dchainfour}[7]{\xymatrix@C-6pt{\overset{#1}{\underset{\ }{\circ}}\ar
@{-}[r]^{#2}
& \overset{#3}{\underset{\ }{\circ}}\ar  @{-}[r]^{#4}  & \overset{#5}{\underset{\
}{\circ}} \ar  @{-}[r]^{#6}
& \overset{#7}{\underset{\ }{\circ}}}}

\newcommand{\Dchainfive}[9]{\xymatrix@C-6pt{\overset{#1}{\underset{\ }{\circ}}\ar
@{-}[r]^{#2}  & \overset{#3}{\underset{\ }{\circ}}\ar  @{-}[r]^{#4}  &
\overset{#5}{\underset{\ }{\circ}}
\ar  @{-}[r]^{#6}  & \overset{#7}{\underset{\ }{\circ}}\ar  @{-}[r]^{#8}  &
\overset{#9}{\underset{\ }{\circ}}}}

\newcommand{\Dtriangle}[6]{
\xymatrix@R-12pt{  &    \overset{#2}{\circ} \ar  @{-}[dl]_{#4}\ar  @{-}[dr]^{#5} & \\
\overset{#1}{\circ} \ar  @{-}[rr]^{#6}  &  &\overset{#3}{\circ} }}

\newcommand{\Dthreefork}[8]{
\rule[-9\unitlength]{0pt}{12\unitlength}
\begin{picture}(28,12)(0,9)
\put(2,10){\ifthenelse{\equal{#1}{l}}{\circle*{2}}{\circle{2}}}
\put(3,10){\line(1,0){10}}
\put(14,10){\ifthenelse{\equal{#1}{m}}{\circle*{2}}{\circle{2}}}
\put(15,10){\line(1,1){7}}
\put(15,10){\line(1,-1){7}}
\put(22,18){\ifthenelse{\equal{#1}{t}}{\circle*{2}}{\circle{2}}}
\put(22,2){\ifthenelse{\equal{#1}{b}}{\circle*{2}}{\circle{2}}}
\put(2,12){\makebox[0pt]{\scriptsize #2}}
\put(8,11){\makebox[0pt]{\scriptsize #3}}
\put(14,12){\makebox[0pt]{\scriptsize #4}}
\put(19,16){\makebox[0pt][r]{\scriptsize #5}}
\put(19,4){\makebox[0pt][r]{\scriptsize #6}}
\put(24,17){\makebox[0pt][l]{\scriptsize #7}}
\put(24,2){\makebox[0pt][l]{\scriptsize #8}}
\end{picture}}

\newcommand{\Drightofway}[8]{\xymatrix@R-6pt{  &    \overset{#6}{\circ} \ar
@{-}[d]_{#4}\ar  @{-}[dr]^{#7} & \\
\overset{#1}{\circ} \ar  @{-}[r]^{#2}  &\overset{#3}{\circ} \ar  @{-}[r]^{#5}
&\overset{#8}{\circ} }}

\usepackage{color}

\numberwithin{equation}{section}\theoremstyle{plain}

\newtheorem{stepo}{Step}

\newtheorem{stepooo}{Step}

\newtheorem{theorem}{Theorem}[section]

\newtheorem{lemma}[theorem]{Lemma}

\newtheorem{conjecture}[theorem]{Conjecture}
\newtheorem{pro}[theorem]{Proposition}

\newtheorem*{lemma*}{Lemma}

\newtheorem{claim}{Claim}

\newtheorem{prop}[theorem]{Proposition}

\theoremstyle{definition}

\newtheorem{definition}[theorem]{Definition}

\newtheorem{example}[theorem]{Example}

\newtheorem{question}[theorem]{Question}

\theoremstyle{remark}
\newtheorem{rem}[theorem]{Remark}
\newtheorem{remark}[theorem]{Remark}

\newcommand{\ydh}{{}^{H}_{H}\mathcal{YD}}
\newcommand{\ydzt}{{}^{\Z^{\theta}}_{\Z^{\theta}}\mathcal{YD}}
\newcommand{\ydg}{{}^{\ku\Gamma}_{\ku\Gamma}\mathcal{YD}}

\newcommand\Ss{\mathcal{S}}

\newcommand\zt{\Z^{\theta}}
\newcommand\G{\mathbb{G}}

\newcommand\Gp{\mathbb{G}'}

\newcommand{\I}{\mathbb{I}}

\newcommand\id{\operatorname{id}}

\newcommand\ord{\operatorname{ord}}
\newcommand\Hom{\operatorname{Hom}}

\newcommand\GK{\operatorname{GK-dim}}

\newcommand\End{\operatorname{End}}

\newcommand\rk{\operatorname{rank}}

\newcommand\car{\operatorname{char}}

\newcommand\gr{\operatorname{gr}}
\newcommand\co{\operatorname{co}}

\newcommand\ad{\operatorname{ad}}

\newcommand\ku{\Bbbk}

\newcommand{\qsr}{\mathfrak{K}}
\newcommand\bq{\mathfrak{q}}
\newcommand{\wbq}{\widetilde{\bq}}
\newcommand\bp{\mathfrak{p}}

\newcommand{\ngo}{\mathfrak{n}}

\newcommand\ot{\otimes}

\newcommand\Z{\mathbb{Z}}
\newcommand\N{\mathbb{N}}

\newcommand\mS{\mathcal{S}}

\newcommand\mP{\mathcal{P}}

\newcommand\cA{\mathcal{A}}

\newcommand\bS{\mathfrak{S}}

\newcommand\bD{\mathfrak{D}}

\newcommand\g{\mathfrak{g}}
\newcommand\h{\mathfrak{h}}

\newcommand{\cE}{{\mathcal E}}

\newcommand{\pre}{\mathfrak{Pre}}
\newcommand{\pref}{\mathfrak{Pre}_{\textrm{fGK}}}
\newcommand{\post}{\mathfrak{Post}}
\newcommand{\postf}{\mathfrak{Post}_{\textrm{fGK}}}

\newcommand{\mbfa}{{}_{\cA}\mathbf{f}}
\newcommand{\mbfk}{{}_{\ku}\mathbf{f}}
\newcommand{\mbfkt}{{}_{\ku}\widetilde{\mathbf{f}}}
\newcommand{\mbfkto}{{}_{\ku_0}\widetilde{\mathbf{f}}}
\newcommand{\mbf}{\mathbf{f}}

\newcommand\pf{\begin{proof}}
\newcommand\epf{\end{proof}}

\begin{document}
\title[Finite GK-dimensional pre-Nichols algebras]
{Finite GK-dimensional pre-Nichols algebras of quantum linear spaces and of  Cartan type}

\author[Nicol\'as Andruskiewitsch and Guillermo Sanmarco]
{Nicol\'as Andruskiewitsch and Guillermo Sanmarco}

\thanks{This material is based upon work supported by the National Science Foundation under Grant No. DMS-1440140 while N. A.
was in residence at the Mathematical Sciences Research Institute in Berkeley, California, in the Spring 2020 semester.
The work of N. A. and G. S. was partially supported by CONICET and Secyt (UNC)}

\address{ Facultad de Matem\'atica, Astronom\'ia y F\'isica,
Universidad Nacional de C\'ordoba. CIEM -- CONICET. 
Medina Allende s/n (5000) Ciudad Universitaria, C\'ordoba, Argentina}
\email{andrus|gsanmarco@famaf.unc.edu.ar}

\begin{abstract}
We study pre-Nichols algebras of quantum linear spaces and of Cartan type with finite GK-dimension. 
We prove that out of a short list of  exceptions involving only roots of order 2, 3, 4, 6, any such pre-Nichols algebra is a quotient of the distinguished pre-Nichols algebra introduced by Angiono generalizing the De Concini-Procesi quantum groups.
There are two new examples, one of which  can be thought of as $G_2$ at a third root of one.
\end{abstract} 
\maketitle

\setcounter{tocdepth}{1}
\tableofcontents

\section{Introduction}
\subsection{Overview}
Let $\ku$ be a field. Let $\GK$ be an abbreviation of Gelfand-Kirillov dimension, see \cite{KL}. 
In this paper we contribute to the ongoing program of classifying Hopf algebras with finite $\GK$. 
See \cite{BGLZ,G-survey,Liu} and references therein.

Let $H$ be a Hopf algebra and let $\ydh$ be the category of Yetter-Drinfeld  modules over $H$.  
Assume that $H$ is pointed (similar arguments apply more generally if its coradical is a Hopf subalgebra).
Basic  invariants of $H$  are 
\begin{enumerate}[leftmargin=*,label=\rm{(\roman*)}]
\item the group of grouplikes $\Gamma = G(H)$,
\item the diagram $R = \oplus_{n \in \N_0} R^n$,  a graded connected Hopf algebra in $\ydg$,
\item the infinitesimal braiding  $V :=  R^1$,  an object in $\ydg$.
\end{enumerate}
See \cite{AS1}. Assume that $\Gamma$ has finite growth.
In order to classify those $H$ with finite $\GK$, one first needs to understand all such $R$ 
with finite $\GK$. Since $R$ is coradically graded and connected it is strictly graded as in \cite{sweedler}.
Strictly graded Hopf algebras $R$ in $\ydg$ with $R^1 \simeq V$ are called post-Nichols algebras of $V$; 
also, graded Hopf algebras $R$ in $\ydg$ generated by $R^1 \simeq V$ are called pre-Nichols algebras of $V$.
See \S \ref{subsec:preliminaries-pre-post-nichols}.

\medbreak
The Nichols algebra $\toba(V)$ is isomorphic to the subalgebra of $R$ generated by $V$.
When  $\car\ku =0$  and $\dim H< \infty$, it was conjectured in \cite{AS-Advances} that 
$R = \toba(V)$, which reduces our problem to classifying finite-dimensional Nichols algebras in $\ydg$.
The conjecture was proved to be valid in \cite{An-diagonal} assuming that $\Gamma$ is abelian. But beyond those hypotheses
this fails to be true.  Thus, it does not seem to be avoidable to consider the following questions:
\begin{enumerate}[leftmargin=*,label=\rm{(\Alph*)}]
\item\label{item-intro-question-A} classify all $V \in \ydg$ such that $\toba(V)$ has finite $\GK$,
\item\label{item-intro-question-B} for such $V$ classify all  post-Nichols algebras with finite $\GK$.
\newcounter{nameOfYourChoice}
\setcounter{nameOfYourChoice}{\value{enumi}}
\end{enumerate}

Lemma \ref{lemma:prenichols-finite-gkd-postnichols} below from \cite{AAH-jordan}
reduces Question \ref{item-intro-question-B} for $V$ as in  \ref{item-intro-question-A} to

\begin{enumerate}[leftmargin=*,label=\rm{(\Alph*)}]
\setcounter{enumi}{\value{nameOfYourChoice}}
\item\label{item-intro-question-C} 
classify all  pre-Nichols algebras of $V^*$ with finite $\GK$.
\end{enumerate}

As usual it is more flexible to deal with these questions considering classes of braided vector spaces 
rather than  classes of groups $\Gamma$ and correspondingly pre-Nichols algebras as braided Hopf algebras.
For Question \ref{item-intro-question-C} we point out that all pre-Nichols algebras of $V$  form a poset
$\pre(V)$ with $T(V)$ minimal and $\toba(V)$ maximal and those with finite $\GK$ form a saturated subposet $\pref(V)$.
In the extreme case when $\car \ku =0$ and the braiding is the 
usual flip, the Nichols algebra is just the symmetric algebra and the pre-Nichols algebras with finite $\GK$ are
the universal enveloping algebras of the finite-dimensional $\N$-graded Lie algebras generated in degree one. Thus
$\pref(V)$ is hardly computable when $\dim V \geq 2$. Similar considerations are valid when the braiding is the super flip of a super vector space, 
see \S \ref{subsubsec:pre-supersym}.
But if $\dim V = 1$, then  $\pref(V) = \pre(V)$ has obviously a minimal element. 
We say that a pre-Nichols algebra is \emph{eminent} if it is a minimum in $\pref(V)$. See Definition \ref{def:eminent}. 
We shall show that many other Nichols algebras with finite $\GK$ have eminent pre-Nichols algebras. 

\medbreak
From now on we assume that $\ku$ is algebraically closed and $\car \ku =0$. 
In this paper we are concerned with  Question \ref{item-intro-question-C} 
for braided vector spaces $V$ of diagonal type, i.e. with braiding
determined by a matrix $\bq = (q_{ij})_{i, j \in \I}$ with entries in $\ku^{\times}$  where  $\theta \in \N$ and  $\I = \{1, \dots, \theta\}$. 
See \S \ref{subsec:preliminaries-nichols-diag} for precisions.

\medbreak
 First we need to discuss Question \ref{item-intro-question-A} for this class. Finite-dimensional Nichols algebras of
diagonal type, i.e. those with $\GK = 0$, were classified in \cite{H-class} through the notion of (generalized) root system.
More generally the list of all Nichols algebras of diagonal type with finite root system is given in \emph{loc. cit.}
It was conjectured in \cite{AAH-triang} that Nichols algebras of
diagonal type with finite $\GK$ are those with finite root system.
This conjecture was verified in various cases \cite{R quantum groups,AA,AAH-diag}. 
We shall assume in a few proofs that the conjecture is true.

\medbreak
Let $\toba(V)$ be a finite-dimensional Nichols algebra of diagonal type with connected Dynkin diagram.
The distinguished pre-Nichols algebra $\wtoba(V)$  was introduced in \cite{An-diagonal} as a tool for determining the defining relations of $\toba(V)$.
Several aspects of these algebras were established in \cite{An-distinguished}. 
Particularly it was asked in \cite{An-distinguished}  (with the terminology just introduced)
whether or not the distinguished are eminent.

\subsection{The main results}
In the present paper we focus on braided vector spaces of diagonal type of two kinds.
Fix $V$ of diagonal type, with braiding given by the matrix $\bq = (q_{ij})_{i, j \in \I}$.

\subsubsection{Quantum linear spaces}\label{subsubsec:intro-qls}
Here  we assume that $\bq$ satisfies 
$q_{ij} q_{ji} = 1$ for all $i \neq j \in \I$.
The \emph{distinguished pre-Nichols algebra} $\wtoba_{\bq}$  is presented  by generators
$(x_i)_{i \in \I}$ and relations $x_i x_j - q_{ij} x_j x_i$, for all $i \neq j \in \I$. 
We need some notation to state our first Theorem. Set
\begin{align}\label{eq:intro-qls-subsets}
\begin{aligned}
\I^{\infty} &= \{i \in \I \colon q_{ii} \notin \G_{\infty} \},&  \I^N &= \{ i\in \I \colon \ord q_{ii}  = N\},\ N \geq 1,
\\
\I^t &= \bigcup_{N > 3} \I^N, &
\I^{\pm} &= \{i \in \I \colon q_{ii} = \pm 1 \} = \I^1 \sqcup \I^2.  
\end{aligned}\end{align}
Thus $\I = \I^{\pm} \sqcup \I^3 \sqcup \I^t \sqcup \I^{\infty}$. For  $\star\in  \N \cup \{\pm, t,\infty \}$, let 
$V^{\star}$ be the subspace of $V$ spanned by $(x_i)_{i \in \I^{\star}}$ and $\bq^{\star}$ the restriction of $\bq$ to $V^{\star}$. Then 
\begin{align*}
V &=   V^{\pm} \oplus V^{3} \oplus V^{t} \oplus V^{\infty}.  
\end{align*}

The pre-Nichols algebras of $V^{\pm}$ with finite $\GK$ are described in \S \ref{subsubsec:pre-supersym}.

\begin{theorem}\label{thm:qls}
\begin{enumerate}[leftmargin=*,label=\rm{(\alph*)}]
\item \label{item:qls-infty-t-3} The distinguished pre-Nichols algebra $\wtoba (V^{\star})$
is eminent, for  $\star\in \{3, t,\infty \}$.

\item \label{item:qls-decomp} Let $\toba$ be a finite GK-dimensional pre-Nichols algebra of $V$; let 
$\toba^{\pm, 3}$, respectively $\toba^t$, $\toba^\infty$ be the subalgebra of $\toba$ generated by
$V^{\pm} \oplus V^{3}$, respectively $V^{t}$, $V^{\infty}$. Then there is a decomposition
\begin{align}\label{eq:qls-decomp}
\toba &\simeq \toba^{\pm, 3} \underline{\otimes} \toba^t \underline{\otimes} \toba^\infty.
\end{align}

\item \label{item:qls-last} Assume that $V$ has a basis $\{x_1, x_2\}$ with $x_1\in V^3$, $x_2\in V^1$.
Then 
\begin{align*}
\tobaqsr(V) = T(V)/ \langle (\ad_c x_{1})^4(x_2), \, (\ad_c x_{2})^2(x_1) \rangle
\end{align*}
 is an eminent pre-Nichols algebra of $V$  and has $\GK = 6$.
\end{enumerate}
\end{theorem}

Parts \ref{item:qls-infty-t-3} and \ref{item:qls-decomp} follow from Proposition \ref{prop:xij=0-preNichols}.
Part \ref{item:qls-last} is Proposition \ref{prop:distinguished-G2-N=3}. Although $\tobaqsr(V)$ of part \ref{item:qls-last}
 is not the distinguished pre-Nichols algebra
of the quantum plane $V$, it can be thought of as the distinguished one of the braided vector space of Cartan type $G_2$, but degenerated in the sense that the parameter is a primitive third root of unity.
Via suitable bosonizations, $\tobaqsr(V)$ provides new examples of pointed Hopf algebras with finite $\GK$.

By \eqref{eq:qls-decomp} it remains to understand $\toba^{\pm, 3}$ for $\toba\in \pref(V)$.  By Proposition \ref{prop:xij=0-preNichols},
$ \toba^{2, 3} \simeq \toba^2 \underline{\otimes} \toba^3$. 
Towards $ \toba^{1, 3}$ we just know Part \ref{item:qls-last}, see also \S \ref{subsec:more-reductions}.
The next step would be the following:

\begin{question}\label{question:intro-qls133}
Assume that $V = V^1 \oplus V^3$, $\dim V^1 = 1$ and $\dim V^3 = 2$.  Is the distinguished pre-Nichols algebra $\wtoba (V)$
eminent?
\end{question}

\subsubsection{Connected Cartan type}\label{subsubsec:intro-cartan}
Here  $\bq$ is of  finite Cartan type, i.e. 
\begin{align*}
q_{ij}q_{ji}  &= q_{ii}^{a_{ij}} , &  i &\neq j \in \I,
\end{align*}
where $\ba = (a_{ij})_{i,j \in \I}$ is a Cartan matrix of finite type with connected Dynkin diagram. In \S \ref{sec:cartan} we recall the possibilities for such $\bq$. They depend on an root of unity $q$, whose order is denoted by $N$.

\begin{theorem}\label{thm:main-intro}
\begin{enumerate}[leftmargin=*,label=\rm{(\alph*)}]
\item \label{item:Cartan-distinguished-eminent}
The distinguished pre-Nichols algebra $\wtoba_{\bq}$
is eminent except in the following cases: $A_2$ with $N=3$, 
\begin{align}\label{eq:exceptions-mainthm}
\begin{aligned}
A_{\theta},\quad \theta \geq 2,\quad N=2; 
& \quad & D_{\theta}, \quad\theta \geq 4, \quad N=2; 
& \quad & G_2, \quad N=4,6.
\end{aligned}
\end{align}
\item \label{item:Cartan-A2-N=3-eminent} Suppose $\bq$ is of type $A_2$ with $N=3$. Then 
\begin{align*}
\htoba=\ku \langle  x_1, x_2 | x_{1112}, x_{2221}, x_{2112}, x_{1221}\rangle
\end{align*}
is an eminent pre-Nichols algebra of $\bq$, and $\GK \htoba = 5$.
\end{enumerate}
\end{theorem}

This answers (partially) a question in \cite{An-distinguished}. 
The  graded duals of the distinguished pre-Nichols algebras have been presented by generators and relations in \cite{AAR-divpow}.

The proof of \ref{item:Cartan-distinguished-eminent} is given in Lemmas \ref{lem:eminent-B2},  \ref{lem:preNichols-G2=image-distinguished}, 
\ref{lem:eminent-B3-C3}, \ref{lem:distinguished-simplylaced-eminent}, \ref{lem:eminent-F4}, \ref{lem:eminent-Etheta}.
For the cases listed in \eqref{eq:exceptions-mainthm} the determination of the poset $\pref(\bq)$ remains an open problem.
See Section \ref{sec:open} for partial results; answers to
Questions \ref{question:A-2-N=2}, \ref{question:A3-N2}, \ref{question:A-4-N=2},
\ref{question:D4-N2}  and \ref{question:D>4-N2} would shed light  on the issue.

The proof of \ref{item:Cartan-A2-N=3-eminent} is given in Proposition \ref{pro:A-2-N=3-eminent}. The eminent pre-Nichols algebra $\htoba$ is introduced and studied in \S \ref{subsubsec:Cartan-A2-N=3}. There we show that $\htoba$ properly covers the distinguished pre-Nichols algebra $\wtoba_{\bq}$, which has $\GK\wtoba_{\bq} =3$. 

In De Concini-Procesi quantum groups at roots of unity and, more generally, in Angiono's distinguished pre-Nichols algebras, 
powers of root vectors generate a skew-central subalgebra \cite[\S 4.1]{An-distinguished}. 
Our $\htoba$ has a  slightly bigger skew central subalgebra $\cZ$ that fits in an extension 
$\ku \to \cZ \hookrightarrow \htoba \twoheadrightarrow \toba_{\bq} \to \ku$ 
of braided Hopf algebras.
We think that the same arguments may apply to produce new examples 
in the case $D_{\theta}$ at $-1$. 
This will be treated in a sequel.

\section{Preliminaries}
\subsection{Conventions}
For $n \leq m \in \N_0$, put $\I_{n,m} = \{k \in \N_0 \colon n \leq k \leq m\}$ and $\I_m = \I_{1,m}$.
Given a positive integer $N$, 
we denote by $\G_N$ the group of $N$-th roots of unity in $\ku^{\times}$, and by $\Gp_N \subset \G_N$ the subset of those of order $N$. 
The group of all roots of unity is denoted by $\G_{\infty}$ and $\Gp_{\infty} := \G_{\infty} - \{1\}$.

The subalgebra generated by a subset $X$ of an associative algebra is denoted by $\ku \langle X \rangle$.

All Hopf algebras are assumed to have bijective antipode. 
If $H$ is a Hopf algebra, the group of group-like elements is denoted by $G(H)$, while $\mP(H)$ is the subspace of primitive elements.
By $\gr H$ we mean the graded coalgebra associated to the coradical filtration.

If $A$ and $B$ are algebras in $\ydh$, we denote by $A \underline{\otimes} B=(A \otimes B, \mu_{A \underline{\otimes} B})$ the algebra with multiplication $\mu_{A \underline{\otimes} B} = (\mu_A \otimes \mu_B)(\id_A \otimes c_{B,A} \otimes \id_B)$, where 
$\mu_A$ and $\mu_B$ are the multiplications of $A$ and $B$, respectively.

\subsection{Gelfand-Kirillov dimension}
We refer to \cite{KL} for general information on this topic.
The following useful statement is immediate from the definition of $\GK$.
Let $R$ be a ring and let $M = \oplus_{n \in \N_0} M_n$ be a graded $R$-module such that 
each $M_n$ is free of finite rank (we say $M$ is locally finite).
The Poincar\'e series of $M$ is  $P_M = \sum_{n \in \N_0} \rk M_n\, X^n \in \Z[[X]]$.

\begin{lemma}\label{lema:poincare-GK} Let $\mathbb L$ and $\mathbb F$ be fields and 
let $T = \oplus_{n \in \N_0} T_n$ and $U = \oplus_{n \in \N_0} U_n$ be two locally finite 
graded algebras generated in degree one over $\mathbb L$ and $\mathbb F$ respectively.
If $P_T  = P_U$, then $\GK T = \GK U$. \qed
\end{lemma}

Actually \cite[12.6.2]{KL} shows that the Poincar\'e series of a graded finitely generated algebra provides its $\GK$.

\subsection{Braided Hopf algebras}

A  pair $(V,c)$ where $V$ is a vector space and $c \in GL(V^{\ot 2})$ satisfies 
the braid equation
\begin{align*}
(c\ot \id) (\id \ot c) (c\ot \id) = (\id \ot c) (c\ot \id) (\id \ot c)
\end{align*}
is called a braided vector space.
A braided vector space with compatible algebra and coalgebra structures as in \cite{Tak}
is called a braided Hopf algebra. 
For instance the tensor algebra $T(V)$ has a canonical structure of (graded connected)  braided Hopf algebra 
such that the elements of degree 1 are primitive. 
Also the tensor coalgebra $T^c(V)$ becomes a braided Hopf algebra by the twisted shuffle product.
There is a homogeneous morphism of braided Hopf algebras  $\Omega \colon T(V) \to T^c(V)$ 
determined by $\Omega (v) = v$, $v \in V$; 
its image is the Nichols algebra  of $V$, denoted $\toba(V)$.

Another description: let $\cJ(V)$
be the largest element of the set $\bS$ of graded Hopf ideals of $T(V)$ trivially intersecting $\ku \oplus V$. 
Then $\toba(V) \simeq T(V)/\cJ(V)$.

\subsection{Principal realizations}
Theorems \ref{thm:qls} and  \ref{thm:main-intro} are relevant for the classification of Hopf algebras with finite $\GK$. 
Indeed a braided vector space arises (up to a mild condition) as a Yetter-Drinfeld module over a Hopf algebra; 
this is called a realization. Realizations are not unique and we single out a class of them for braidings of diagonal type.
Let $H$ be a Hopf algebra. 
A YD-pair is a couple
$(g, \chi) \subset G(H) \times \Hom_{\text{Alg}} (H, \ku)$ satisfying 
\begin{align*}
\chi(h) g &= \chi(h_{(2)}) h_{(1)} g \mS(h_{(3)}), && h \in H.
\end{align*}
Compare with \cite[p. 671]{AS-p3}.
This compatibility guarantees that $\ku^{\chi}_g$ (i.~e. $H$ acting and coacting on $\ku$ by $\chi$ and $g$, respectively) is a Yetter-Drinfeld
module over $H$.
Let $(V,c^{\bq})$ be a braided vector space of diagonal type. Following \cite[p. 673]{AS-p3}, a principal realization of $(V,c^{\bq})$ over $H$
is a family $(g_i, \chi_i )_{ i \in \I}$ of YD-pairs such that $q_{ij} = \chi_j(g_i)$ for all $i, j$. In this case $V = \bigoplus_i \ku^{\chi_i}_{g_i} \in \ydg$.

\subsection{Pre-Nichols and post-Nichols algebras}\label{subsec:preliminaries-pre-post-nichols}
We present in detail the objects of interest in this paper.

\begin{itemize} [leftmargin=*]
\item Let $\toba =  \bigoplus_{n \in \N_0} \toba^n$ be a graded connected braided Hopf algebra 
with $\toba^1 \simeq V$. Then $\toba$ is a \emph{pre-Nichols} algebra of $V$ 
if it is generated by $\toba^1$.
In this case there are epimorphisms of (graded) braided Hopf algebras 
\begin{align*}
T(V) \twoheadrightarrow \toba \twoheadrightarrow \toba(V).
\end{align*}
\end{itemize}
Hence the set $\pre(V)$ of isomorphism classes of pre-Nichols algebras of $V$ is partially ordered 
with $T(V)$ minimal and $\toba(V)$ maximal:
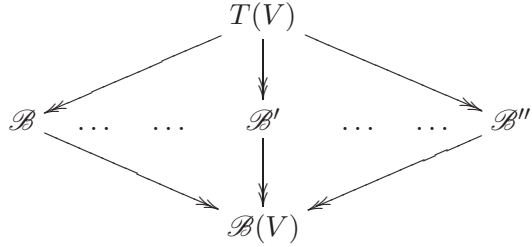
\begin{align*}
\xymatrix@C-15pt{ &&& T(V) \ar @{->>}[1,-3] \ar @{->>}[1,0] \ar @{->>}[1,3]  &&& 
\\ \toba & \dots & \dots &\toba' & \dots & \dots & \toba'' 
\\ & & & \toba(V) \ar @{<<-}[-1,-3] \ar @{<<-}[-1,0] \ar @{<<-}[-1,3]  & & &}
\end{align*}

\begin{itemize} [leftmargin=*]
\item Dually, a graded connected braided Hopf algebra $\cE = \bigoplus_{n \in \N_0} \cE^n$  with $\cE^1 \simeq V$ is a \emph{post-Nichols} algebra of $V$ if it is coradically graded.
Thus we have monomorphisms of (graded) braided Hopf algebras 
\begin{align*}
\toba(V) \hookrightarrow \cE \hookrightarrow T^c(V).
\end{align*}
\end{itemize}

Hence the set $\post(V)$ of isomorphism classes of post-Nichols algebras of $V$ is partially ordered 
with $T^c(V)$ maximal  and $\toba(V)$ minimal:
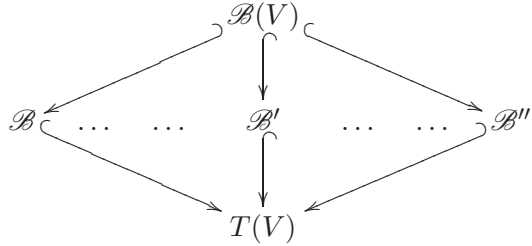
\begin{align*}
\xymatrix@C-15pt{ &&& \toba(V) \ar @{_{(}->}[1,-3] \ar @{^{(}->}[1,0] \ar @{^{(}->}[1,3]  &&& 
\\ \toba & \dots & \dots &\toba' & \dots & \dots & \toba'' 
\\ & & & T(V) \ar @{<-_{)}}[-1,-3] \ar @{<-_{)}}[-1,0] \ar @{<-^{)}}[-1,3]  & & &}
\end{align*}
The only pre-Nichols which is also a post-Nichols algebra of $V$ is $\toba(V)$ itself.

\subsection{Eminent pre- and  post-Nichols algebras}\label{subsec:eminent}
For the purposes of classifying  Hopf algebras with finite $\GK$, it is important to describe the 
(partially ordered) subset $\postf(V)$ of $\post(V)$ consisting of post-Nichols algebras with  finite $\GK$.
In this paper we are mainly interested in the (partially ordered) subset $\pref(V)$ of $\pre(V)$ consisting of pre-Nichols algebras with  finite $\GK$.
The reason to start with this is given by the following result:

\begin{lemma}\label{lemma:prenichols-finite-gkd-postnichols} \cite{AAH-jordan}
Let $\toba$ be a pre-Nichols algebra of $V$ and let $\cE = \toba^d$ be the graded dual of $\toba$.
Then  $\GK \cE \leq \GK \toba$. If $\cE$ is finitely generated, then the equality holds. \qed
\end{lemma}

A first approximation to the determination of $\postf(V)$ and $\pref(V)$ is through the following notion.

\begin{definition}\label{def:eminent}
\begin{enumerate}[leftmargin=*,label=\rm{(\alph*)}]
\item A pre-Nichols algebra $\htoba$ is \emph{eminent} if it is the minimum of $\pref(V)$; i.~e.  
there is an epimorphism of braided Hopf algebras $\htoba \twoheadrightarrow \toba$ that is the identity on $V$
for any $\toba \in \pref(V)$.

\smallbreak
\item A post-Nichols algebra $\widehat{\cE}$ is \emph{eminent} if it is the maximum of $\postf(V)$; that is for any $\cE \in \postf(V)$, 
there is a monomorphism of braided Hopf algebras $\cE \hookrightarrow \widehat{\cE}$ that is the identity on $V$.
\end{enumerate}
\end{definition}

Beware that there are braided vector spaces without eminent pre-Nichols algebras; e.~g., if $\dim V > 1$ and the braiding is the usual flip, then
$\pref(V)$ has infinite chains. An intermediate situation could be described as follows.

\begin{definition}\label{def:eminent-system}
A family $(\htoba_i)_{i\in I} \subset \pref(V)$  is \emph{eminent} if 
\begin{enumerate}[leftmargin=*,label=\rm{(\alph*)}]
\item \label{item:eminent-system-1} for any $\toba \in \pref(V)$, 
there exists $i\in I$ and  an epimorphism of braided Hopf algebras $\htoba_i \twoheadrightarrow \toba$ that is the identity on $V$, and
\item \label{item:eminent-system-2} $(\htoba_i)_{i\in I}$ is minimal among the families in $\pref(V)$ satisfying \ref{item:eminent-system-1}.
\end{enumerate}
Eminent families of post-Nichols algebras are defined similarly.
\end{definition}

All the notions above about braided Hopf algebras related to braided vector spaces have a counterpart 
for Yetter-Drinfeld modules. Namely, suppose that $(V,c)$ is realized in $\ydh$ for some Hopf algebra $H$. Then
$\pre^H(V)$ is the subset of $\pre(V)$ of pre-Nichols algebras that belong to $\ydh$; similarly
we have  $\pref^H(V)$, $\post^H(V)$, $\postf^H(V)$, and also $H$-eminent pre-Nichols or post-Nichols algebras.

%
%
%
%

\subsection{The adjoint representation and $q$-brackets}
Any Hopf algebra $R$ in $\ydh$ comes equipped with the (left) \emph{adjoint} representation $\ad_c \colon R \to \End R$, given by
\begin{align*}
(\ad_c x) y&= \mu (\mu \ot \Ss) (\id \ot c) (\Delta \ot \id) (x \ot y) , && x, y \in R,
\end{align*}
where $\mu$, $\Delta$ and $\Ss$ denote the multiplication, comultiplication and antipode of $R$, respectively.
The adjoint action of a primitive element $x \in R$ is 
\begin{align*}
(\ad_cx)y&=xy - (x_{(-1)}\cdot y)x_{(0)}, &   y \in R.
\end{align*}
Given $x_{i_1}, x_{i_2} \dots, x_{i_k} \in R$, put
\begin{align}\label{eq:xij-def}
x_{i_1 i_2 \dots i_k} =  (\ad_c x_{i_1}) \dots (\ad_c x_{i_{k-1}}) x_{i_k}.
\end{align}
We also set
$x_{(k \, h)} = x_{k\,(k+1)\, (k+2) \dots h}$ for $k < h$.

On the other hand, the \emph{braided commutator} is defined by
\begin{align*}
[x, y]_c&=xy - (x_{(-1)}\cdot y)x_{(0)}, &   x,y \in R.
\end{align*}
We refer to \cite[Introduction]{AA-diag} for a more detailed treatment.

\subsection{Nichols algebras of diagonal type}\label{subsec:preliminaries-nichols-diag}
Fix a natural number $\theta$ and let $\I = \I_{\theta}$. Any matrix $\bq = (q_{ij})_{i, j \in \I}$ with coefficients in $\ku^{\times}$ 
determines a braided vector space \emph{of diagonal type} $(V, c^{\bq})$, where 
\begin{align}\label{eq:def-diagonal}
V \text{ has a basis } (x_i) _ {i \in \I}, && c^{\bq}(x_i \ot x_j) = q_{ij} x_j \ot x_i, && i, j \in \I. 
\end{align}

The \emph{Dynkin diagram} associated to $\bq$ is a non-oriented graph with $\theta$ vertices. The vertex $i$ is labelled by $q_{ii}$, 
and there is an edge between $i$ and $j$ if and only if $\widetilde{q}_{ij}:= q_{ij}q_{ji} \neq 1$; in this case, the edge is labeled by $\widetilde{q}_{ij}$. 
Thus we may speak of the connected components of this diagram and by abuse of notation of $\bq$.
The following useful result says that a connected component with at least 2 vertices one of them labelled by 1 gives rise to 
an infinite GK-dimensional Nichols algebra.

\begin{lemma} \label{lemma:2dim,qii=1} \cite[Lemma 2.8] {AAH-triang}
Let $U$ be a braided vector space of diagonal type with Dynkin diagram 
\begin{align*}
\xymatrix{ \underset{ \ }{\overset{q}{\circ}} \ar  @{-}[rr]^{r}  & & \underset{ \ }{\overset{ 1 }{\circ} }} , \quad r\neq 1.
\end{align*} 
Then $\GK \toba(U) = \infty$. \qed
\end{lemma}

\medbreak
Let $\alpha_1, \dots, \alpha_{\theta}$ be the canonical basis of $\Z^{\theta}$.
From the braiding matrix $\bq$ we obtain a $\ku^{\times}$-valued bilinear form on $\Z^{\theta}$, still denoted $\bq$ and determined by 
$\bq (\alpha_i, \alpha_j) = q_{ij}$, $i, j \in \I$. Put also
\begin{align} \label{eq:def-qtilde}
\wbq (\alpha, \beta)& := \bq(\alpha, \beta) \bq(\beta, \alpha), & \alpha, \, \beta \in \Z^{\theta}.
\end{align}
For sake of brevity, we use $\bq_{\alpha \beta}=\bq(\alpha, \beta)$ and $\wbq_{\alpha \beta}=\wbq(\alpha, \beta)$ as well.

\medbreak The braided vector space $(V,c^{\bq})$ as in \ref{eq:def-diagonal} is realized in $\ydzt$ by declaring
\begin{align}\label{eq:V-ydzt}
\deg (x_i) &= \alpha_i,& \alpha_i \cdot x_j &= q_{ij}x_j, & i,j &\in \I.
\end{align}
The algebra $T(V)$ becomes $\Z^{\theta}$-graded. 
Thus any quotient  algebra $R$ of $T(V)$ by a graded ideal inherits the grading:
$R = \bigoplus_{\alpha \in \Z^{\theta}} R^{\alpha}$. We keep the notation $\deg$ for this degree.
Furthermore, if $R$ is an algebra obtained as a quotient of $T(V)$ by a graded ideal $I$ (thus a subobject in $\ydzt$), 
then the braiding on the homogeneous subspaces is given by
\begin{align}\label{eq:braid-quotient-tensor}
c(u \ot v) &= \bq_{\alpha, \beta} \, v \ot u, && u\in R^{\alpha}, \, v \in R^{\beta}.
\end{align}
The braided commutators satisfy
\begin{align}
\label{eq:braided-commutator-right-mult}
[u,vw]_c &= [u,v]_c w + \bq_{\alpha \beta} v [u,w]_c,
\\
\label{eq:braided-commutator-left-mult} 
[uv,w]_c &= \bq_{ \beta \gamma}[u,w]_c v + u[v,w]_c,
\\
\label{eq:braided-commutator-iteration}
\big[ [u,v]_c, w \big]_c &= \big[ u, [v, w]_c \big]_c 
- \bq_{\alpha \beta} v [u,w]_c 
+ \bq_{ \beta \gamma}[u,w]_c v,
\end{align}
for homogeneous elements $u\in R^{\alpha}$,  $v \in R^{\beta}$, $w \in R^{\gamma}$.
\medbreak
In the diagonal setting \eqref{eq:def-diagonal} we set as usual $\cJ_{\bq} = \cJ(V)$, $\toba_{\bq} =\toba(V)$, etc. 
Nichols algebras \emph{of diagonal type} (i.~e. those arising from braided vector spaces of diagonal type)
have been intensively studied. The classification of all matrices $\bq$ such that
$\toba_{\bq}$ has finite root system was provided in \cite{H-class}; the defining relations of these Nichols algebras are given 
in \cite{An-diagonal,An-convex}.
Clearly, finite dimensional Nichols algebras of diagonal type have finite root system. 
It was  conjectured that those of finite $\GK$ share the same property.

\begin{conjecture}\label{conj:AAH} \cite [Conjecture 1.5] {AAH-triang} The root system 
of a Nichols algebra of  diagonal type with finite GK-dimension is finite. 
\end{conjecture}

The validity of Conjecture \ref{conj:AAH}  would imply the classification of finite GK-dimensional Nichols algebras of diagonal type. 
There is strong evidence supporting it.
The conjecture holds when $\theta = 2$ \cite[Thm. 4.1] {AAH-diag},  when the braiding is of affine Cartan type \cite[Thm. 1.2] {AAH-diag}, or when
$\bq$ is generic, that is $q_{ii} \notin \G_{\infty}$, and $q_{ij}q_{ji}=1$ or $q_{ij}q_{ji}\notin \G_{\infty}$,
for all $i\ne j\in \I$ \cite{R quantum groups, AA}.

We include for completeness proofs of the following well-known results.

\begin{lemma} \label{lem:disconnected-primitives}
Let $0\neq v, w \in T(V)$ be homogeneous primitive elements with $\deg v = \alpha$ and $\deg w = \beta$.
Then $(\ad_cv)w$ is primitive if and only if $\wbq_{\alpha \beta} = 1$.
\end{lemma}

\noindent \emph{Proof.} 
Using \eqref{eq:braid-quotient-tensor}, compute $\Delta((\ad_cv)w)= \Delta (vw - \bq_{\alpha \beta}  \, wv)=$
\begin{align*}
&= (v \ot 1 + 1 \ot v)(w \ot 1 + 1 \ot w) 
- \bq_{\alpha \beta} \, (w \ot 1 + 1 \ot w)(v \ot 1 + 1 \ot v)\\
&= vw \ot 1 + v\ot w + \bq_{\alpha \beta} \, w\ot v + 1 \ot vw\\
&\quad - \bq_{\alpha \beta} \, (wv \ot 1 + w\ot v + \bq_{\beta \alpha} \, v\ot w + 1 \ot wv)\\
&=(\ad_cv)w \ot 1 + 1 \ot (\ad_cv)w + \big(1-\wbq_{\alpha \beta} \big) v\ot w. \hspace{92pt}\qed
\end{align*}

\begin{lemma} \label{lem:subspace-of-primitives}
Let $R$ be a graded braided Hopf algebra. If $W$ is any braided subspace of $R$ contained in $\mP(R)$ then $\GK \toba(W) \leq \GK R$.
\end{lemma}

\pf
We follow \cite[Lemma 5.4]{AS-Annals}.
Since the elements of $W$ are primitive, the subalgebra $\ku \langle W \rangle$ is a braided Hopf subalgebra of $R$; by definition of the Nichols algebra it follows that $\gr \ku \langle W \rangle$ projects onto $\toba(W)$, so $\GK \toba(W) \leq \gr \ku \langle W \rangle$. But $\GK \gr \ku \langle W \rangle \leq \GK \ku \langle W \rangle$ by \cite[Lemma 6.5]{KL}, and this proves the desired inequality.
\epf

\subsection{Pre-Nichols algebras of diagonal type}\label{subsec:preliminaries-prenichols-diag}
Let $(V, c^{\bq})$ be a braided vector space of diagonal type associated to the matrix
$\bq = (q_{ij})_{i, j \in \I}$. Recall that $\widetilde{q}_{ij} = q_{ij}q_{ji}$, $i\neq j$.
We write $\pref^{\zt} (V)$ for $\pref^{\ku \zt} (V)$, cf. \eqref{eq:V-ydzt}.

\subsubsection{Pre-Nichols algebras under twist-equivalence}
Let $\bp = (p_{ij})_{i, j \in \I}$ be another braiding matrix such that
\begin{align*}
q_{ii} &= p_{ii},& \widetilde{q}_{ij}  &= \widetilde{p}_{ij},& i, j &\in \I.
\end{align*}
In this case, $(V, c^{\bq})$ and the braided vector space $(W, c^{\bp})$ with basis $(y_i)_{i\in \I}$
are said to be {\it twist-equivalent}.

\begin{lemma}\label{lema:pref-twisting}
There is an isomorphism of posets $\pref^{\zt} (W) \simeq \pref^{\zt} (V)$.
\end{lemma}

\pf  Let $\sigma: \zt \times \zt \to \ku^{\times}$ be the  bilinear form, hence a 2-cocycle, given by
$\sigma(\alpha_{i}, \alpha_{j}) = \begin{cases} p_{ij}  q^{-1}_{ij}, & i\le j, \\
1, & i > j\end{cases}$.
Let $T(V)_{\sigma}$ be the corresponding cocycle deformation of $T(V)$, i.~e. with multiplication
\begin{align}\label{eq:twistingalgebra}
u._{\sigma}v &= \sigma(\alpha, \beta) u v,&  u&\in T(V)^{\alpha}, \, v \in T(V)^{\beta}, \, \alpha, \beta \in \zt.
\end{align}
By the proof of \cite[Prop. 3.9]{AS1} the linear map $\varphi: W \to V$, 
$\varphi(y_{i}) = x_{i}$, $i\in \I$, induces an isomorphism   $\varphi: T(W) \to T(V)_{\sigma}$
of Hopf algebras in $\ydzt$.  Let $I$ be a Hopf ideal of $T(V)$ that belongs to $\ydzt$; then 
it is also a Hopf ideal of $T(V)_{\sigma}$ and $\GK T(V)/ I = \GK T(V)_{\sigma}/ I$ 
by Lemma \ref{lema:poincare-GK}.
\epf

\subsubsection{Pre-Nichols algebras of super symmetric algebras}\label{subsubsec:pre-supersym}
Assume that $\widetilde{q}_{ij} = 1 = q_{ii}^2$, for all $j\neq i \in \I$.
Then $V = V_0 \oplus V_1$ is a super vector space where
$V_{j}$ is spanned by those $x_i$'s such that $q_{ii} = (-1)^j$, $j =0, 1$.
Let  $\bp = (p_{ij})_{i, j \in \I}$ be the matrix corresponding to the associated super symmetry.
Then 
\begin{itemize} [leftmargin=*]
\item The pre-Nichols algebras of $(V, c^{\bp})$ are  the enveloping superalgebras $U(\ngo)$, where $\ngo = \oplus_{j\in \N} \ngo^j$ is a graded Lie superalgebra generated by $\ngo^1 \simeq V$. 

\item  $\pref (V, c^{\bp})$ consists of  the enveloping superalgebras $U(\ngo)$, 
where $\ngo = \oplus_{j\in \N} \ngo^j$ is a graded Lie superalgebra generated by $\ngo^1 \simeq V$ with $\dim \ngo < \infty$.

\item Hence $\pref^{\zt} (V, c^{\bp})$ consists of  the enveloping superalgebras $U(\ngo)$, 
where $\ngo  = \oplus_{\beta\in \zt} \ngo^\beta$ is a finite-dimensional $\zt$-graded Lie superalgebra 
generated by $\ngo^1 = \oplus_{i\in \I} \ngo^{\alpha_{i}} \simeq V$. In particular $\pref^{\zt} (V, c^{\bp}) \subsetneq \pref (V, c^{\bp})$.

\item By Lemma \ref{lema:pref-twisting}, $\pref^{\zt} (V, c^{\bq})$ is isomorphic as a poset to the 
set of isomorphism classes of finite-dimensional $\zt$-graded Lie superalgebras as in the previous point.
\end{itemize}

\subsubsection{Distinguished pre-Nichols algebras}
Assume that $\dim \toba_{\bq} < \infty$.
The \emph{distinguished} pre-Nichols algebra of $V$ introduced in \cite{An-distinguished} 
is the quotient $\wtoba_{\bq} := T(V)/ \cI_{\bq}$, where $\cI_{\bq}$ is the ideal of $T(V)$ generated by the defining relations of $\cJ_{\bq}$ given in \cite{An-diagonal} but excluding the powers of root vectors and including the quantum Serre relations at Cartan vertices.
A detailed presentation of $\cJ_{\bq}$ and $\cI_{\bq}$ is available in \cite[\S 4]{AA-diag}.

\section{Quantum linear spaces}\label{sec:qls}
In this section we investigate finite GK-dimensional pre-Nichols algebras of quantum linear spaces. These are Nichols algebras of braided vector spaces of diagonal type with totally disconnected Dynkin diagram.
More precisely, fix a matrix $\bq = (q_{ij})_{i,j \in \I}$ and a vector space $V$ with basis $(x_i)_{i \in \I}$  and braiding given by $c^{\bq}(x_i \ot x_j) = q_{ij} x_j \ot x_i$, $i, j \in \I$. In this section we assume that 
\begin{align} \label{eq:def-QLS}
q_{ij} q_{ji} = 1, && i \neq j \in \I.
\end{align} 
Then $\tobaq$ is presented by generators
$(x_i)_{i \in \I}$ and relations 
\begin{align} \label{eq:B(V)-rels}
x_{ij} &= 0, \, \, \text{ if } i<j,  & & \\ 
x_i^{N_i} &= 0, \, \, \text{ if } q_{ii} \in \Gp_{\infty},& \text{where } N_i &:= \ord q_{ii} \in \N \cup \infty;
\end{align}
here we are using the notation \eqref{eq:xij-def}. It has a PBW-basis:
\begin{align} \label{eq:B(V)-PBW}
\{ x_1^{ a_1} x_2^{ a_2} \cdots x_\theta ^ { a_{\theta}} \colon 0 \leq a_i < N_i \text { if } q_{ii} \in \Gp_{\infty}; \, \,   0\leq a_i \text{ otherwise} \} .
\end{align}

\medbreak 
The \emph{distinguished pre-Nichols algebra} $\wtoba_{\bq}$ of $V$ is presented  by generators
$(x_i)_{i \in \I}$ and relations \eqref{eq:B(V)-rels}; it is a domain of $\GK = \theta$. Recall the partition  
$\I = \I^{\pm} \sqcup \I^3 \sqcup \I^t \sqcup \I^{\infty}$ where as in \eqref{eq:intro-qls-subsets} we set
\begin{align*}
\I^{\infty} &= \{i \in \I \colon q_{ii} \notin \G_{\infty} \},&  
\I^{\pm} &= \{i \in \I \colon q_{ii} = \pm 1 \}, \\
\I^N &= \{ i\in \I \colon q_{ii} \in \Gp_N\},\ N \geq 3,& 
\I^t &= \bigcup_{N > 3} \I^N.  
\end{align*}
 For  $\star\in \{\pm, 3, t,\infty \}$, let 
$V^{\star}$ be the subspace of $V$ spanned by $(x_i)_{i \in \I^{\star}}$ and $\bq^{\star}$ the restriction of $\bq$ to $V^{\star}$. Then 
$V =   V^{\pm} \oplus V^{3} \oplus V^{t} \oplus V^{\infty}$.  
As we have seen in \S \ref{subsubsec:pre-supersym} the $\zt$-graded pre-Nichols algebras of $V^{\pm}$ are twistings of 
enveloping algebras of nilpotent Lie superalgebras with suitable properties, particularly there is no eminent 
pre-Nichols algebra of $V^{\pm}$.


\subsection{Reduction to order $\leq 3$}\label{subsec:reductions-3}

\begin{remark} \label{rmk:W-Dynkin}
Let $i \neq j \in \I$. Recall that $x_{ij} := (\ad_c x_i) x_j= x_i x_j - q_{ij} x_j x_i$. The braiding of the 3-dimensional subspace $\ku x_i+\ku x_{ij} + \ku x_j  \subset T(V)$ is easily computed, and the corresponding Dynkin diagram is either
\begin{align}\label{diagram:QLS-i-ij-j}
\xymatrix{ \underset{ i }{\overset{q_{ii}}{\circ}} \ar  @{-}[rr]^{q_{ii}^2}  & & \underset{ ij }{\overset{ q_{ii} q_{jj} }{\circ} } 
\ar  @{-}[rr]^{q_{jj}^2}  & & \underset{ j }{\overset{ q_{jj} }{\circ} }},
\end{align}
or it is disconnected if the label of some edge is $1$.
\end{remark}

\begin{pro} \label{prop:xij=0-preNichols}
Let $i, j \in \I$ such that $4< \ord q_{ii} + \ord q_{jj}$.
Then $x_{ij} = 0$ holds in any finite GK-dimensional pre-Nichols algebra of $\bq$. 
\end{pro}

\pf
Let $\toba $ be a pre-Nichols algebra of $\bq$, so there is a braided Hopf algebra map $T(V) \to \toba$. 
Let $y_1, y_2, y_3$ denote the image of $x_i, x_j, x_{ij}$, respectively, and consider $W := \ku y_1 + \ku y_2 + \ku y_3$. 
By Lemma \ref{lem:disconnected-primitives} we have $W\subset \mP (\toba)$, hence Lemma \ref{lem:subspace-of-primitives} warranties $\GK \toba(W) \leq \GK \toba$.

Assume $y_3 \neq 0$, so $W$ is 3-dimensional by a degree argument and its Dynkin diagram $\bD$ is \eqref{diagram:QLS-i-ij-j}. We show that $\GK\toba(W) = \infty$. 

Consider the subspaces $ V_1 = \ku y_1 \oplus \ku y_3$, $ V_2 = \ku y_3 \oplus \ku y_2 \subset W$; 
denote their corresponding Dynkin diagrams by $\bD _1 $ and $\bD _2$, respectively. 
From $q_{ij} q_ {ji} = 1$ it follows $x_{ij} = - q_{ij} x_{ji}$, so the image of $x_{ji}$ in $\toba$ is not zero.

We split the proof in several cases according to the possibilities for $\ord q_{ii}$ and $\ord q_{jj}$.

\medbreak
\textbf{Case 1: $q_{ii} \notin \G_{\infty}$ or $q_{jj} \notin \G_{\infty}$}. This essentially goes back to \cite{R quantum groups}.
Assume first $q_{ii} \notin \G_{\infty}$.
If $\GK \toba < \infty$, it follows from \cite [Lemmas 2.6 and 2.7] {AAH-triang} 
that there exists a natural number $k$ such that $(k)^! _ {q_{ii}} \prod_{h=0}^{k-1} (1- q^h_{ii}) = 0$, which contradicts $q_{ii} \notin \G_{\infty}$. 
The case $q_{jj} \notin \G_{\infty}$ is  similar: since the image of $x_{ji}$ is not zero, 
we may apply the same argument as with $q_{ii}$. 

\textbf{Case 2: $q_{jj} = 1$}. By the previous case, we might assume $q_{ii}$ is a root of unity, and by hypothesis its order must be $N_i > 3$.
The  diagram $\bD_1 $ is 
\begin{align*}
\xymatrix{ \underset{ i }{\overset{q_{ii}}{\circ}} \ar  @{-}[rr]^{q_{ii}^2}  & & \underset{ ij }{\overset{q_{ii} }{\circ} }}, &&\text{Cartan type }\quad 
\begin{pmatrix}
2 & 2 - N_i  \\
2 - N_i & 2
\end{pmatrix}.
\end{align*}
If $\GK \toba(V_1) < \infty$ then \cite{AAH-diag} implies that the Cartan matrix is of finite type. 
Thus we conclude $N_i=3$, a contradiction.

\textbf{Case 3: $q_{jj} = -1$}. By Case 1, assume that $q_{ii}$ is a root of unity. Its order is $\geq3$. 
By  \cite {AAH-diag}, $\GK \toba (V_1) = \infty$ since the Dynkin diagram of 
$V _1$ is 
$$ \xymatrix{ \underset{ i }{\overset{q_{ii}}{\circ}} \ar  @{-}[rr]^{q_{ii}^2}  & & \underset{ ij }{\overset{- q_{ii} }{\circ} }} $$
and this does not appear in \cite [Table 1] {H-class}.

\textbf{Case 4: $q_{ii},q_{jj} \notin \Gp_2$}. Now $W$ has a connected  Dynkin diagram:
$$ \bD = \xymatrix{ \underset{ i }{\overset{q_{ii}}{\circ}} \ar  @{-}[rr]^{q_{ii}^2}  & & \underset{ ij }{\overset{ q_{ii} q_{jj} }{\circ} } 
\ar  @{-}[rr]^{q_{jj}^2}  & & \underset{ j }{\overset{ q_{jj} }{\circ} }} .$$
If the Nichols algebra of $V_1$ is finite GK-dimensional, by exhaustion of \cite [Table 1] {H-class} we conclude that $q_{ii}$, $q_{jj}$ and $\bD_1$ satisfy one of the following:
\begin{multicols}{2}
\begin{enumerate} [leftmargin=*] 
\item \label{diag:rank2row2}
{\small $q_{ii} \in \Gp_3, q_{ii}q_{jj} = -1$, }
$\xymatrix{ \underset{ i }{\overset{q_{ii}}{\circ}} \ar  @{-}[r]^{q_{ii}^2}  &  \underset{ ij }{\overset{-1}{\circ} }}$
\item \label{diag:rank2row3a}
$ q_{ii} \in \Gp_4$, $q_{jj} = q_{ii},
\xymatrix{ \underset{ i }{\overset{q_{ii}}{\circ}} \ar  @{-}[r]^{-1}  &  \underset{ij}{\overset{-1}{\circ}}}$
\item \label{diag:rank2row3b}
{\small $ q_{ii} \in \Gp_3, q_{jj} = \pm q_{ii}$, }
$\xymatrix{ \underset{i}{\overset{q_{ii}}{\circ}} \ar  @{-}[r]^{q_{ii} ^ 2}  & \underset{ij}{\overset{\pm q_{ii}^2 }{\circ} }} $
\item \label{diag:rank2row5}
$q_{ii} \in \Gp_3, q_{jj} =  q_{ii}^2, 
\xymatrix{ \underset{i}{\overset{q_{ii}}{\circ}} \ar  @{-}[r]^{q_{ii}^2}  & \underset{ij}{\overset{ 1 }{\circ} }} $
\item \label{diag:rank2row7}
{\small $q_{ii} \in \Gp_6,  q_{ii} q_{jj} = -1$, }
$\xymatrix{ \underset{i}{\overset{q_{ii}}{\circ}} \ar  @{-}[r]^{q_{ii}^2}  & \underset{ij}{\overset{ -1 }{\circ} }} $
\item \label{diag:rank2row10a}
$ q_{ii} \in \Gp_5, q_{jj} = q_{ii} ^ 2, \xymatrix{ \underset{i}{\overset{q_{ii}}{\circ}} \ar  @{-}[r]^{q_{ii}^2}  & \underset{ij}{\overset{ q_{ii} ^ 3 }{\circ} }} $
\item \label{diag:rank2row10b}
$ q_{jj}  \in \Gp_9, q_{ii} = q_{jj} ^ 3, 
\xymatrix{ \underset{i}{\overset{q_{jj} ^ 3}{\circ}} \ar  @{-}[r]^{q_{jj}^6}  & \underset{ij}{\overset{ q_{jj} ^ 4 }{\circ} }} $
\item \label{diag:rank2row13}
{\small $ q_{ii} \in \Gp_5, q_{ii} q_{jj} = -1$, }
$\xymatrix{ \underset{i}{\overset{q_{ii}}{\circ}} \ar  @{-}[r]^{q_{ii}^2}   & \underset{ij}{\overset{ -1 }{\circ} }} $
\end{enumerate}
\end{multicols}

In the rest of the proof, we discard one by one all these possibilities.

\eqref{diag:rank2row3a} Now $W$ is of Cartan type with Dynkin diagram and Cartan matrix:
\begin{align*}
\bD = \xymatrix{ \underset{i}{\overset{q_{ii}}{\circ}} \ar  @{-}[rr]^{-1}  & & \underset{ij}{\overset{ -1 }{\circ} } \ar  @{-}[rr]^{ -1}  & & \underset{j}{\overset{ q_{ii} }{\circ} }}, \,  q_{ii} \in \Gp_4;
&&\begin{pmatrix}
2 &  -2  & 0 \\
-1 &  2  & -1 \\
0 & -2 & 2
\end{pmatrix}.
\end{align*}
Since this matrix is of affine type, $\GK \toba(W)=\infty$ by \cite {AAH-diag}. 

\eqref{diag:rank2row3b} Assume first $q_{jj} = - q_{ii}$. Then
$$\bD_2 = \xymatrix{ \underset{ij}{\overset{ - q_{ii} ^ 2}{\circ}} \ar  @{-}[rr]^{q_{ii} ^ 2}  & & \underset{j}{\overset{- q_{ii} }{\circ} }}, \quad  q_{ii} \in \Gp_3, $$
which is not arithmetic. By \cite {AAH-diag} we see that $\GK \toba(V_2) = \infty$. 
Next, when $q_{jj} =  q_{ii}$, $W$ is of Cartan type with  Dynkin diagram and Cartan matrix:
\begin{align*}
\bD = \xymatrix{ \underset{i}{\overset{q_{ii}}{\circ}} \ar  @{-}[rr]^{q_{ii}^ 2} & & \underset{ij}{\overset{ q_{ii} ^ 2 }{\circ} } 
\ar  @{-}[rr]^{ q_{ii} ^ 2}  & & \underset{j}{\overset{ q_{ii} }{\circ} }},\,   q_{ii} \in \Gp_3; &&\begin{pmatrix} 2 &  -1  & 0 \\-2 &  2  & -2 \\0 & -1 & 2\end{pmatrix},
\end{align*}
which is affine, so $\GK \toba (W) = \infty $ by \cite{AAH-diag}.
\eqref{diag:rank2row5} Since $q_{ii} ^ 2 \neq 1$, we have $\GK \toba(V_1)=\infty$ by \cite[Lemma 2.8] {AAH-triang}.
\eqref{diag:rank2row10a} In this case
$$\bD_2 = \xymatrix{ \underset{ij}{\overset{q_{ii}^3}{\circ}} \ar  @{-}[rr]^{q_{ii}^4}  & & \underset{j}{\overset{ q_{ii} ^ 2 }{\circ} }}, \quad q_{ii} \in \Gp_5, $$
is of indefinite Cartan type, so  $\GK \toba (V_2) = \infty$ by \cite {AAH-diag}.
\eqref{diag:rank2row10b} Similarly, 
$$\bD_2 = \xymatrix{ \underset{ij}{\overset{q_{jj} ^ 4}{\circ}} \ar  @{-}[rr]^{q_{jj}^2}  & & \underset{j}{\overset{ q_{jj}  }{\circ} }},  \quad q_{jj} \in \Gp_9, $$
is indefinite Cartan, so $ \GK \toba (V_2) = \infty$.
In the remaining cases,  $ \mathfrak D$ is 

\eqref{diag:rank2row2}  
$ \bD = \xymatrix{ \underset{i}{\overset{\omega }{\circ}} \ar  @{-}[rr]^{\omega ^ 2}  & & \underset{ij}{\overset{ - 1}{\circ} } 
\ar  @{-}[rr]^{\omega}  & & \underset{j}{\overset{ - \omega ^ 2 }{\circ} }}, \quad \omega \in \Gp_3 .$

\eqref{diag:rank2row7}  
$ \bD = \xymatrix{ \underset{i}{\overset{ - \omega ^ 2 }{\circ}} \ar  @{-}[rr]^{\omega }  & & \underset{ij}{\overset{ - 1}{\circ} } 
\ar  @{-}[rr]^{\omega ^ 2}  & & \underset{j}{\overset{  \omega  }{\circ} }}, \quad \omega \in \Gp_3 .$

\eqref{diag:rank2row13}  
$ \bD = \xymatrix{ \underset{i}{\overset{\zeta }{\circ}} \ar  @{-}[rr]^{\zeta ^ 2}  & & \underset{ij}{\overset{ - 1}{\circ} } 
\ar  @{-}[rr]^{\zeta ^ 3 }  & & \underset{j}{\overset{ - \zeta ^ 4 }{\circ} }}, \quad \zeta \in \Gp_5 .$

Now \eqref{diag:rank2row2} and \eqref{diag:rank2row7}  are equal up to permutation of the indexes.
Only here we need to assume the validity of Conjecture \ref{conj:AAH}. 
Indeed, these diagrams do not appear in \cite [Table 2] {H-class}, so   $\GK \toba (W)= \infty$ in all cases.
\epf

\subsection{A pre-Nichols algebra of type $G_2$}
We assume $(V,c^{\bq})$ has the following Dynkin diagram
\begin{align*}
\xymatrix{ \underset{ 1 }{\overset{\omega}{\circ}}  & & \underset{ 2 }{\overset{ 1 }{\circ} }} , \qquad \omega \in \Gp_3.
\end{align*} 

\begin{prop} \label{prop:distinguished-G2-N=3}
The algebra $\tobaqsr_{\bq} := T(V)/ \langle x_{11112}, \, x_{221} \rangle$ 
is an eminent pre-Nichols algebra of $(V, c^{\bq})$ and $\GK \tobaqsr_{\bq} = 6$.
\end{prop}

\pf
We first claim that the elements $x_{11112}$ and $x_{221}$ are primitive in $T(V)$.
This is verified by a direct computation, see \cite{sanmarco}.

Second, we claim that the relations $x_{11112}=0$ and $x_{221}=0$ hold in any finite GK-dimensional pre-Nichols algebra  $\toba$ of $(V, c^{\bq})$.

Assume first $x_{11112}\neq0$ in $\toba$. Then also $x_{12}\neq0$. From Lemma \ref{lem:disconnected-primitives} and the previous claim, we have a braided subspace 
\begin{gather*}
W=\ku x_1 +\ku x_{12}+\ku x_{11112} \subset \mP(\toba),
\end{gather*}
so Lemma \ref{lem:subspace-of-primitives} gives $\GK \toba(W) \leq \GK \toba$. By a degree argument, $W$ has dimension three; from direct computation its Dynkin diagram is
\begin{align}\label{eq:q-serre-G2-N=3-y11112}
\begin{aligned}
\xymatrix{ &\underset{1}{\overset{\omega}{\circ}} \ar @{-}[rd]^{\omega^2}  \ar @{-}[ld]_{\omega^2} &  \\ \underset{12}{\overset{\omega}{\circ}} \ar  @{-}[rr]^{\omega^2}  & &\underset{11112}{\overset{\omega}{\circ}}, }
\end{aligned} 
& &\text { Cartan type } 
\begin{pmatrix} 2 & -1 & -1\\ -1&2&-1\\-1&-1 & 2 \end{pmatrix}. 
\end{align}
Since the Cartan matrix is of affine type $A^1_2$, we have $\GK \toba(W) = \infty$ by  \cite[Theorem 1.2(a)]{AAH-diag}. Thus $\GK \toba = \infty$. 

Assume now $x_{221}\neq0$ in $\toba$. Then $x_{21}\neq0$, and since $q_{12}q_{21}=1$, we have $x_{12}=-q_{12} x_{21}\neq0$. Consider $W'=\ku x_1 +\ku x_{12}+\ku x_{221}$. We may now use the same argument as above. Indeed, $W' \subset \mP(\toba)$ has Dynkin diagram \eqref{eq:q-serre-G2-N=3-y11112} replacing $x_{11112}$ by $x_{221}$, so $\GK \toba(W') = \infty$ by the same reason as $\GK \toba(W) = \infty$. Hence $\GK \toba = \infty$. 
Thus $\tobaqsr_{\bq} \twoheadrightarrow \toba$.

The verification of $\GK \tobaqsr_{\bq} = 6$ is postponed to Proposition \ref{prop:tobaqsr-q-q^3}.
\epf

\subsection{A further reduction}\label{subsec:more-reductions}
Let $\toba$ be a finite GK-dimensional pre-Nichols algebra of $\toba_{\bq}$. We are naturally led to consider
\begin{align}
E := \{(i,j) \colon \, i \in \I^3, j \in \I^1, \, x_{ij} \neq 0 \text{ in } \toba \}.
\end{align}

\begin{remark} \label{rmk:xi-xij-A2}
If $(i,j) \in E$, the braided vector space $\ku x_i \oplus \ku x_{ij}\subset \toba$ is of Cartan type $A_2$ by Remark \ref{rmk:W-Dynkin}. 
\end{remark}

\begin{lemma}\label{lema:Eij}
If $(i, j_1), (i,j_2) \in E$ then $j_1=j_2$. 
\end{lemma}
\pf
Since $x_{ij_1}$ and $x_{ij_2}$ are $\Z^{\theta}$-homogeneous,
\begin{align*}
c(x_{ij_1}\ot x_{ij_2}) = \bq(\alpha_i + \alpha_{j_1}, \alpha_i + \alpha_{j_2}) \, x_{ij_2}\ot x_{ij_1}= q_{ii}q_{ij_2}q_{j_1i}q_{j_1j_2} \, x_{ij_2}\ot x_{ij_1}.
\end{align*}
Assume $j_1 \neq j_2$. Then $x_i$, $x_{ij_1}$ and $x_{ij_2}$ have pairwise different $\Z^{\theta}$-degrees, 
so they span a 3-dimensional braided subspace $W = \ku x_i  + \ku x_{ij_1} + \ku  x_{ij_2} \subset \mP(\toba)$. 
Now the Dynkin diagram of $W$ is 
\begin{align*}
\begin{aligned}
\xymatrix{ &\underset{i}{\overset{q_{ii}}{\circ}} \ar @{-}[rd]^{q_{ii}^2}  \ar @{-}[ld]_{q_{ii}^2} &  \\ \underset{ij_1}{\overset{q_{ii}}{\circ}} \ar  @{-}[rr]^{q_{ii}^2}  & &\underset{ij_2}{\overset{q_{ii}}{\circ}}, }
\end{aligned} 
&  q_{ii} \in \Gp_3,&\text {Cartan type } 
\begin{pmatrix} 2 & -1 & -1\\ -1&2&-1\\-1&-1 & 2 \end{pmatrix}. 
\end{align*}
Since the Cartan matrix is of affine type $A^{(1)}_2$, we have $\GK \toba(W) = \infty$ by  \cite[Theorem 1.2(a)]{AAH-diag}. Thus $\GK \toba = \infty$, a contradiction. 
\epf

\section{Cartan type}\label{sec:cartan}
In this section we determine the finite GK-dimensional pre-Nichols algebras of braided vector spaces of finite Cartan type under some restrictions.

We fix a matrix $\bq = (q_{ij})_{i,j \in \I}$ of non-zero scalars such that $q_{ii} \neq 1$ for all $i\in \I$ 
and a braided vector space $(V, c^{\bq})$  
with braiding given by $c^{\bq}(x_i \ot x_j) = q_{ij} x_j \ot x_i$, $i, j \in \I_{\theta}$, in 
a basis $\{x_1, \dots, x_{\theta}\}$.  Let $N_i=\ord q_{ii} \in \N \cup \infty$.

\medbreak
Recall that $\bq$, or $(V,c^{\bq})$, is of \emph{Cartan type} if there exists a Cartan matrix $\ba = (a_{ij})_{i,j \in \I}$ such that  $q_{ij}q_{ji}=q_{ii}^{a_{ij}}$ for all $i,j$. 
Let $i \in \I$. If $N_i=\infty$, then $a_{ij}$ are uniquely determined. Otherwise, we impose 
\begin{align}\label{eq:cartan-type}
N_i &< a_{ij} \leq 0,& j & \in \I.
\end{align}
In this way we say that $(V,c^{\bq})$, is of Cartan type $\ba$.

We follow the terminology of \cite{K}. Cartan matrices are arranged in three families, namely: finite, affine and indefinite. We say that $\bq$, or $(V,c^{\bq})$, belongs to one of these families if the corresponding $\ba$ does. 

In this section we assume that $\bq$ is of \emph{connected finite Cartan type}. Here are the possibilities for the Dynkin diagram of $\bq$:
\begin{align*}
&A_{\theta}: \ \xymatrix{ \overset{q}{\underset{\ }{\circ}}\ar  @{-}[r]^{q^{-1}}  &
 \overset{q}{\underset{\
}{\circ}}\ar@{.}[r] & \overset{q}{\underset{\ }{\circ}} \ar  @{-}[r]^{q^{-1}}  &
\overset{q}{\underset{\ }{\circ}}}, &
&B_{\theta}: \ \xymatrix{ \overset{\,\,q^2}{\underset{\ }{\circ}}\ar  @{-}[r]^{q^{-2}}  &
\overset{\,\,q^2}{\underset{\ }{\circ}}\ar@{.}[r] & \overset{\,\,q^2}{\underset{\
}{\circ}} \ar  @{-}[r]^{q^{-2}}  & \overset{q}{\underset{\ }{\circ}}}
\\
&C_{\theta}:  \xymatrix{ \overset{q}{\underset{\ }{\circ}}\ar  @{-}[r]^{q^{-1}}  &
\overset{q}{\underset{\ }{\circ}}\ar  @{-}[r]^{q^{-1}} &  \overset{q}{\underset{\
}{\circ}}\ar@{.}[r] & \overset{q}{\underset{\ }{\circ}} \ar  @{-}[r]^{q^{-2}}  &
\overset{\,\,q^2}{\underset{\ }{\circ}}}, &&
\\
&D_{\theta}:  \begin{aligned}
\xymatrix@R-15pt{ 
& & &   \overset{q}{\circ} &\\
\overset{q}{\underset{\ }{\circ}}\ar  @{-}[r]^{q^{-1}}  & 
\overset{q}{\underset{\ }{\circ}}\ar@{.}[r] &
\overset{q}{\underset{\ }{\circ}} \ar  @{-}[r]^{q^{-1}}  & \overset{q}{\underset{\ }{\circ}} \ar @<0.7ex> @{-}[u]_{q^{-1}}^{} \ar  @{-}[r]^{q^{-1}} &
\overset{q}{\underset{\ }{\circ}}}
\end{aligned},&
\\
& E_{\theta}:  \begin{aligned}
\xymatrix@R-15pt
{ &  &   \overset{q}{\circ} &\\
\overset{q}{\underset{\ }{\circ}}\ar  @{-}[r]^{q^{-1}}  &  \overset{q}{\underset{\ }{\circ}}\ar@{.}[r]  & 
\overset{q}{\underset{\ }{\circ}} \ar @<0.7ex> @{-}[u]_{q^{-1}}^{} \ar
@{-}[r]^{q^{-1}} &  
\overset{q}{\underset{\ }{\circ}}  \ar  @{-}[r]^{q^{-1}} &
\overset{q}{\underset{\ }{\circ}}}
\end{aligned}, && \theta \in \I_{6, 8},
\\
&F_4: \ 
\xymatrix{ 
\overset{\,\,q}{\underset{\ }{\circ}}\ar  @{-}[r]^{q^{-1}}  &
\overset{\,\,q}{\underset{\ }{\circ}}\ar  @{-}[r]^{q^{-2}} &   \overset{q^2}{\underset{\ }{\circ}} \ar  @{-}[r]^{q^{-2}}  &  \overset{q^2}{\underset{\ }{\circ}} },
& &
G_2: \ \xymatrix{  
\overset{\,\,q}{\underset{\ }{\circ}} \ar  @{-}[r]^{q^{-3}}&
\overset{q^3}{\underset{\ }{\circ}}}.
\end{align*}
Here $q$ is a root of unity in $\ku$; set $N=\ord q$. 
We refer to the survey \cite{AA-diag} for restrictions on $N$ and  precise features of $\tobaq$ on each case.

The quantum Serre relations are the following elements of $T(V)$:
\begin{align}\label{eq:def-q-Serre}
(\ad_cx_i)^{1-a_{ij}}x_j, &&i, \,j \in \I, & \, \, i\neq j.
\end{align} 
By \cite[Lemma A.1]{AS-Advances} these are primitive in any pre-Nichols algebra.
Let $\wtoba_{\bq} = T(V)/ \cI_{\bq}$ be the distinguished pre-Nichols algebra of $(V,c^{\bq})$, see 
\S \ref{subsec:preliminaries-prenichols-diag}. 
\begin{rem}\label{rem:qsr-presentation}
From the detailed presentation in \cite[\S 4]{AA-diag} we see 
that the quantum Serre relations \eqref{eq:def-q-Serre} generate $\cI_{\bq}$ in the following cases:
\begin{itemize} 
\item when $\ba$ is of type $A_2$ or $B_2$ \cite[pp. 397, 399, 400]{AA-diag},
\item when $\ba$ is of type $G_2$ and $N\neq4, 6$ \cite[pp. 410, 411]{AA-diag},
\item when $\ba$ is simply-laced and $N > 2$ \cite[pp. 397, 404, 407]{AA-diag},
\item when $\ba$ is of type B, C, or F and $N > 4$ \cite[pp. 399, 402, 409]{AA-diag}.
\end{itemize}
\end{rem}

\subsection{Quantum Serre relations}\label{subsub:QSR}
Let $\ba = (a_{ij})_{i,j \in \I}$ be a symmetrizable indecomposable 
generalized Cartan matrix and $\mathbf{d}\in \operatorname{GL}_{\theta} (\Z)$ diagonal such that
$\mathbf{d} \ba$ is symmetric; this is equivalent to an 
irreducible Cartan datum as in \cite[1.1.1]{L-libro} by setting
\begin{align*}
\cdot: \I &\times \I \to \Z, & i\cdot j &= d_{i}a_{ij}, & i,j &\in \I.
\end{align*}
Let $\g = \g(\ba)$ be the associated Kac-Moody algebra which have a triangular decomposition
$\g(\ba) = \g^{+} \oplus \h \oplus \g^{-}$.

Let $q \in \ku^{\times}$ and consider the Dynkin diagram 
\begin{align}\label{eq:Dinkin-Cartan}
\xymatrix{ & \overset{q^{d_i}}{\underset{i}{\circ}}  \ar  @{..}[l] \ar  @{-}[rr]^{q^{d_{i}a_{ij}}} & &
\overset{q^{d_j}}{\underset{j}{\circ}}  \ar  @{.}[r]    & }
\end{align}
Let $\bq$ be any matrix with Dynkin diagram \eqref{eq:Dinkin-Cartan} and $(V,c^{\bq})$ be the 
corresponding braided vector space with basis $(x_i)_{i\in \I}$.  Notice that $\bq$ is of Cartan type
but it is not necessarily of type $\ba$ as \eqref{eq:cartan-type} may not hold.

Let $\tobaqsr_{\bq} = T(V)$ modulo the ideal 
$\qsr_{\bq}$ generated by the quantum Serre relations $(\ad_c x_i)^{1-a_{ij}} (x_j)$, $i\neq j \in \I$, which is a pre-Nichols algebra of $V$.

\begin{prop}\label{prop:tobaqsr-GK} $\GK \tobaqsr_{\bq} \geq \dim \g^+$.
\end{prop}

\pf If $\xi \in \ku$, $\xi^2 = q$, 
then $\bp = (\xi^{d_ia_{ij}})$ has Dynkin diagram \eqref{eq:Dinkin-Cartan}. 
Let $(W,c^{\bq})$ be the 
corresponding braided vector space with basis $(\widehat{x}_i)_{i\in \I}$.

\begin{claim}\label{step:QST-twisting}  $\GK \tobaqsr_{\bq} = \GK \tobaqsr_{\bp}$. 
\end{claim}

\pf By the proof of \cite[Proposition 3.9]{AS1} (or the proof of Lemma \ref{lema:pref-twisting}) there is a  homogeneous linear isomorphism $\psi: T(V) \to T(W)$
determined by $\psi(x_i) = \widehat{x}_i$ for all $i\in \I$ and 
satisfying \cite[Remarks 3.10]{AS1}. Hence $\psi(\qsr_{\bq}) = \qsr_{\bp}$ and $\psi$ induces a homogeneous
linear isomorphism $\psi: \tobaqsr_{\bq} \to \tobaqsr_{\bp}$. Then apply Lemma \ref{lema:poincare-GK}.
\epf

Let now $\mbf$ be the $\mathbb Q(v)$-algebra defined in \cite[1.2.5]{L-libro}, where 
$v$ is an indeterminate and let $\mbfa$ be the $\cA := \Z[v, v^{-1}]$-subalgebra spanned by the quantum divided powers of the generators of $\mbf$
\cite[1.4.7]{L-libro}. By \cite[14.4.3]{L-libro}, $\mbfa$ is a free $\cA$-module and
\begin{align}\label{eq:poincare-fA}
P_{\mbfa} &= P_{\mbf}.
\end{align}
Consider $\ku$ as $\cA$-module via $v \mapsto \xi$. Then we have the algebras $\mbfk = \ku \otimes_{\cA} \mbfa$
and $\mbfkt$ defined in \cite[33.1.1]{L-libro} (which is nothing else than $\tobaqsr_{\bp}$).
By \cite[1.4.3]{L-libro}, the quantum Serre relations hold in $\mbfk$, hence we have a surjective algebra map
$\tobaqsr_{\bp} = \mbfkt \twoheadrightarrow \mbfk$. Thus 
\begin{align}\label{eq:poincare-qsr}
\GK \tobaqsr_{\bp} &\geq \GK \mbf.
\end{align}
On the other hand, let $\ku_0$ be $\ku$ 
as $\cA$-module via $v \mapsto 1$. Then $\mbfkto \simeq U(\g^{+})$ by \cite[33.1.1]{L-libro} and
$\dim_{\mathbb Q(v)} \mbf_{\nu} = \dim_{\ku_0} (\mbfkto_{\nu})$ by \cite[33.1.3]{L-libro}; that is
\begin{align}\label{eq:poincare-U(g)}
\GK \mbf &=\GK( \mbfkto)= \dim \g^+,
\end{align}
where the first equality holds by Lemma \ref{lema:poincare-GK}. The Proposition follows.
\epf

\begin{example}\label{exa:toba6}
Let $\ba = \begin{pmatrix}
2 & -5 \\ -1 & 2 
\end{pmatrix}$. Then \eqref{eq:Dinkin-Cartan} takes the form $\xymatrix{ \overset{q}{\underset{1}{\circ}}\ar  @{-}[r]^{q^{-5}}  &
\overset{q^5}{\underset{2 }{\circ}}}$ with
$q \in \ku^{\times}$. If $q_{12} \in \ku^{\times}$ and $q_{21} := q_{12}^{-1}q^{-5}$, 
then $ \bq = \begin{pmatrix} q & q_{12} \\ q_{21} &q^{5} \end{pmatrix}$ 
has the  Dynkin diagram above.
Here $\tobaqsr_{\bq} = \ku \langle  x_1, x_2 \rangle$ modulo the relations
\begin{align*}
x_2^2x_1 &-q_{21}(q)_2 \, x_2x_1x_2+q q_{21}^2 \ x_1 x_2^2,\\
x_1^6x_2 &- 3q_{12}^2 \ x_1^4 x_2x_1^2 +3q_{12}^4 \ x_1^2x_2x_1^4 - q_{12}^6 \ x_2x_1^6 .
\end{align*} 
In this setting Proposition \ref{prop:tobaqsr-GK} gives $\GK \tobaqsr_{\bq} = \infty$.
\end{example}

\begin{example}\label{exa:tobaG2}
Let $\ba = \begin{pmatrix}2 & -3 \\ -1 & 2 \end{pmatrix}$. 
Then \eqref{eq:Dinkin-Cartan} takes the form 
$\xymatrix{ \overset{q}{\underset{1}{\circ}}\ar  @{-}[r]^{q^{-3}}  &
\overset{q^3}{\underset{2 }{\circ}}}$ with
$q \in \ku^{\times}$. If $q_{12} \in \ku^{\times}$ and $q_{21} := q_{12}^{-1}q^{-3}$, 
then $ \bq = \begin{pmatrix} q & q_{12} \\ q_{21} &q^{3} \end{pmatrix}$ 
has the  Dynkin diagram above.
Here $\tobaqsr_{\bq} = \ku \langle  x_1, x_2 \rangle$ modulo the relations
\begin{align*}
&x_2^2x_1 -q_{21}(2)_{q^3} \, x_2x_1x_2+q_{21}^2 q^3 \ x_1 x_2^2,\\ 
&x_1^4 x_2 -q_{12} (4)_q \, x_1^3x_2x_1 + q_{12}^2 q  \binom{4}{2}_{q}  \, x_1^2x_2x_1^2 -q_{12} ^3 (4)_q \, x_1^3x_2x_1 + q_{12}^4 q^6 \, x_2 x_1^4.
\end{align*} 
In this situation Proposition \ref{prop:tobaqsr-GK} establishes $\GK \tobaqsr_{\bq} \geq 6$.
\end{example}

This last example gains more relevance when the parameter $q\in\ku^{\times}$ specializes to a root of unity with small order.

\begin{prop} \label{prop:tobaqsr-q-q^3} 
Let $\ba$ and $q\in\ku^{\times}$ as in Example \ref{exa:tobaG2}.
\begin{enumerate}[leftmargin=*,label=\rm{(\alph*)}]
\item \label{item:eminent-G2-N=3} If $q\in\Gp_3$ then $\GK \tobaqsr_{\bq} = 6$. 
\item  \label{item:eminent-A2-N=2} If $q\in\Gp_2$ then $x_{112}^2=0$ in $\tobaqsr_{\bq}$. 
\end{enumerate}
\end{prop}
\pf
Put $x_{11122} = [x_{112},x_{12}]$. By direct computation, in $\tobaqsr_{\bq}$ the following relations hold:
\begin{align*}
x_{12} x_{2} &= q_{12} q^3 x_{2} x_{12},   \\
x_{112}x_{2}&= q_{12}^2 q^3x_{2}x_{112}+ q_{12} q^2(q-1)(2)_q x_{12}^2,\\
x_{1112}x_{2}&= \text{\small $q_{12}^3 q ^3$} x_{2}x_{1112}+ \text{\small $q_{12} q(q^2-q-1)$} x_{11122}  + \text{\small $q_{12}^2 q^2(q-1)(3)_q$} x_{12}x_{112},\\
x_{11122}x_{2}&=q_{12}^3 q^6 x_{2}x_{11122} +q_{12}^2q^3(q-1)^2(2)_q x_{12}^3,\\
x_{1}x_{11122}&= q_{12}^2 q^3 x_{11122}x_{1} + x_{1112}x_{12} - q_{12}^2 q^3 x_{12}x_{1112},\\
(2)_qx_{1112}x_{12}&= q_{12}^2 q^3(2)_qx_{12}x_{1112}+  q_{12}q(q-1)(3)_qx_{112}^2.
\end{align*}
\medbreak
\ref{item:eminent-G2-N=3} In this case the last two relations above become
\begin{align*}
x_{1}x_{11122}= q_{12}^2  x_{11122}x_{1},  &&
x_{1112}x_{12}= q_{12}^2 x_{12}x_{1112}.
\end{align*}
These imply more commutations:
\begin{align*}
x_{11122} x_{12} = q_{12} x_{12} x_{11122},&&
x_{112} x_{11122} &= q_{12} x_{11122} x_{112},\\
x_{1112} x_{112} = q_{12} x_{112} x_{1112}, &&
x_{1112}x_{11122}&= q_{12}^3 x_{11122}x_{1112}.
\end{align*}

\smallbreak
Now we claim that $\tobaqsr_{\bq}$ is linearly spanned by  
\begin{align*}
\mathbb{B}=\{x_2^{n_1} x_{12}^{n_2} x_{11122}^{n_3}x_{112}^{n_4}x_{1112}^{n_5}x_1^{n_6}\colon 0\leq n_1, \dots, n_6 \}.
\end{align*}
Denote by $\cI$ the linear span of $\mathbb{B}$. Since $1\in \cI$, it is enough to show that $\cI$ is left ideal of $\tobaqsr_{\bq}$. If we multiply $x_2^{n_1} x_{12}^{n_2} x_{11122}^{n_3}x_{112}^{n_4}x_{1112}^{n_5}x_1^{n_6}$ by $x_1$ on the left, we can use the previously deduced commutations between the (powers of the) $x_{\alpha}$'s to successively rearrange the terms until we get a linear combination of elements in $\mathbb{B}$. 
The claim follows.

\smallbreak
Consider the $\mathbb{N}_0$-filtration $\mathcal{F}$ on $\tobaqsr_{\bq}$ induced by $\mathbb{B}$, and denote by $\gr_{\mathcal F} (\tobaqsr_{\bq})$ the associated graded algebra. There is a natural projection from a polynomial algebra $\ku [y_1, \dots, y_6] \twoheadrightarrow \gr_{\mathcal F} (\tobaqsr_{\bq})$, hence $\GK \gr_{\mathcal F} (\tobaqsr_{\bq})\leq 6$. By \cite[Proposition 6.6]{KL} and Example \ref{exa:tobaG2} we also have $\GK \gr_{\mathcal F} (\tobaqsr_{\bq}) = \GK \tobaqsr_{\bq} \geq 6$, so the equality holds.

\medbreak
\ref{item:eminent-A2-N=2} This follows by specialization at $q=-1$ in the relation \\
$ \qquad \qquad
(2)_qx_{1112}x_{12}= q_{12}^2 q^3(2)_qx_{12}x_{1112}+  q_{12}q(q-1)(3)_qx_{112}^2.
$
\epf

\begin{remark}\label{rem:A2-N=2-missing-PBWgenerator}
Let us point out the relevance of \ref{item:eminent-A2-N=2}. By Kharchenko's theory \cite{Kh}, $\tobaqsr_{\bq}$ has a PBW-basis. 
By Proposition \ref{prop:tobaqsr-GK} we know $\GK \tobaqsr_{\bq} \geq 6$ but, when $q\in\Gp_2$, the root $2\alpha_1+\alpha_2$ will not contribute to $\GK \tobaqsr_{\bq}$ by  \ref{item:eminent-A2-N=2} above.  So even if $\ba$ is of type $G_2$, one should not expect that the PBW generators are just those related to the six positive roots of $G_2$, as was the case in the proof of \ref{item:eminent-G2-N=3}. 
\end{remark}

\subsection{Type $A_2$} \label{subsec:Cartan-A2}
In this and the next subsections we seek for eminent (families of) pre-Nichols algebras in order to determine finite $\GK$ pre-Nichols algebras of braidings of finite Cartan type. The distinguished pre-Nichols algebra will serve as the principal guide in our exploration.

\subsubsection{Type $A_2$ with $N>3$} \label{subsubsec:Cartan-A2-N>3}

\begin{lemma}\label{lem:preNichols-A2-N>3-eminent}
Assume $\ba$ is of Cartan type $A_2$ with $N>3$. If $\toba$ is a finite GK-dimensional pre-Nichols algebra of $\bq$, then $x_{112}=0$ and $x_{221}=0$ in $\toba$, i.~e. the distinguished pre-Nichols algebra $\wtoba_{\bq}$ is eminent, cf. Definition \ref{def:eminent}.
\end{lemma}

\pf Assume $x_{iij}\neq0$ for some $i\neq j \in \I_2$; the 3-dimensional braided subspace  $W:=\ku x_j \oplus \ku x_i \oplus \ku x_{iij} \subset \mP(\toba)$ has $\GK \toba(W) < \infty$. 

\smallbreak
Consider the braided subspace $W_1 = \ku x_i \oplus \ku x_{iij} \subset W$. By direct computation, the braiding on $W_1$ is of Cartan type with the following Dynkin diagram and Cartan matrix: 
\begin{align*}
\xymatrix{ \underset{ i }{\overset{q_{ii}}{\circ}} \ar  @{-}[rr]^{q_{ii}^3}  & & \underset{ iij }{\overset{q_{ii}^3 }{\circ} }}, &&
A_1=\begin{pmatrix} 2 & 3-N\\ 1-M & 2 \end{pmatrix}, && 
M=\begin{cases}
N/3, &\text { if } 3 | N,\\
N, &\text { otherwise}.
\end{cases}
\end{align*}
If either $N=5$ or $N>6$, it is evident that the Cartan matrix $A_1$ is not finite, so $\GK \toba(W_1)=\infty$ by \cite[Theorem 1.2 (b)]{AAH-diag}. This contradicts $\GK \toba < \infty$. 

\noindent
For the remaining cases (i.~e. $N=4$ and $N=6$), we consider the whole $W$. Since 
$ \wbq(\alpha_j, 2\alpha_i+\alpha_j)=(q_{ji}q_{ij})^2q_{jj}^2=1$,
the braiding on $W$ is of Cartan type with the following Dynkin diagram and Cartan matrix
\begin{align*}
\xymatrix{ \underset{j}{\overset{q_{ii}}{\circ}} \ar  @{-}[rr]^{q_{ii}^{-1}}  & & \underset{i}{\overset{ q_{ii} }{\circ} } \ar  @{-}[rr]^{q_{ii}^3}  & & \underset{iij}{\overset{ q_{ii}^3 }{\circ} }},
&&A= \begin{pmatrix} 2 & -1 & 0\\ -1&2&3-N\\0&1-M & 2 \end{pmatrix}.
\end{align*}
Now it is straightforward to verify that if $N=4$ or $6$, then $A$ is of affine type, which contradicts $\GK \toba(W) < \infty$ by \cite[Theorem 1.2 (b)]{AAH-diag}.
\epf

\subsubsection{Type $A_2$ with $N=3$} \label{subsubsec:Cartan-A2-N=3}
Here is the first restriction.

\begin{lemma}\label{lem:A2-N=3-discard}
Assume $\ba$ is of Cartan type $A_2$ with $N=3$. Let $\toba\in\pref$. Then $x_{iiij}=0$ and $x_{jiij}=0$ in $\toba$ for all $i\neq j \in \I_2$. 
\end{lemma}

\pf Since $x_{iij}$ is primitive, using that $\wbq(\alpha_i, 2\alpha_i+\alpha_j)=q^4q^{-1}=1$ and
$\wbq(\alpha_j, 2\alpha_i+\alpha_j)=q^2q^{-2}=1$,
we get $x_{iiij}, x_{jiij} \in \mP(\toba)$ by Lemma \ref{lem:disconnected-primitives}. Assume first $x_{iiij}\neq0$ in $\toba$. The braided subspace $\ku x_i \oplus \ku x_j \oplus \ku x_{iiij} \subset \mP(\toba)$ has finite $\GK$ Nichols algebra. The Dynkin diagram is 
\begin{align*}
\begin{aligned}
\xymatrix{ &\underset{iiij}{\overset{q}{\circ}} \ar @{-}[rd]^{q^{-1}}  \ar @{-}[ld]_{q^{-1}} &  \\ \underset{i}{\overset{q}{\circ}} \ar  @{-}[rr]^{q^{-1}}  & &\underset{j}{\overset{q}{\circ}}, }
\end{aligned} 
& &\text {Cartan type } 
\begin{pmatrix} 2 & -1 & -1\\ -1&2&-1\\-1&-1 & 2 \end{pmatrix}. 
\end{align*}
The Cartan matrix is of affine type $A^{(1)}_2$, and by \cite[Theorem 1.2(a)]{AAH-diag} this contradicts $\GK \toba < \infty$.

\smallbreak If $x_{jiij}\neq0$, the same argument leads to a contradiction. Indeed, by direct computations, the Dynkin diagram of $U=\ku x_i \oplus \ku x_j \oplus \ku x_{jiij}$ is
\begin{align*}
\xymatrix{ &\underset{jiij}{\overset{q}{\circ}} \ar @{-}[rd]^{q^{-1}}  \ar @{-}[ld]_{q^{-1}} &  \\ \underset{i}{\overset{q}{\circ}} \ar  @{-}[rr]^{q^{-1}}  & &\underset{j}{\overset{q}{\circ}}. }
\end{align*}
so $\GK \toba(U)=\infty$ by \cite[Theorem 1.2(a)]{AAH-diag}.
\epf


\begin{remark}\label{rmk:A-2-N=3-eminent}
Denote $\htoba=T(V)/\langle   x_{1112}, x_{2221}, x_{2112}, x_{1221}\rangle$. 
The defining ideal of $\htoba$ is a Hopf ideal by the proof of Lemma \ref{lem:A2-N=3-discard}. 
Let $\pi\colon \htoba \to \toba(V)$ denote the natural projection.  Let $\cZ$ be the subalgebra of $\htoba$ generated by
\begin{align*}
z_1 := x_2^3, \qquad z_2 := x_{221}, 
\qquad z_3 := x_{112}, \qquad z_4 := x_1^3, \qquad z_5 := x_{12}^3.
\end{align*}
The next results are devoted to prove that $\htoba$ is eminent. 
\end{remark}

\begin{lemma}\label{lem:A-2-N=3-central}
\begin{enumerate}[leftmargin=*,label=\rm{(\alph*)}]
\item \label{item:A-2-N=3-eminent-computations} Given $i\neq j \in \I_2$, the following relations hold in $\htoba$:
\begin{align*}
[ x_{ij}, x_{iij}] = 0 = [ x_{ij}, x_{jji}]; \qquad [ x_{iij}, x_{jji}]=0.    
\end{align*}
\item \label{item:A-2-N=3-central-normal} 
$\cZ$ is a normal Hopf subalgebra of $\htoba$.
\item \label{item:A-2-N=3-central-basis} 
The $z_i$'s $q$-commute; $B= \left\{ z_1^{n_1} z_2^{n_2} z_3^{n_3} z_4^{n_4} z_5^{n_5} \colon n_i \in \N_0 \right\}$
is a basis of $\cZ$.
\item \label{item:A-2-N=3-central-coinvariants}
$\cZ = {}^{\co \pi }\htoba$.
\end{enumerate}
\end{lemma}

\pf
\ref{item:A-2-N=3-eminent-computations} Just compute using \eqref{eq:braided-commutator-iteration}:
\begin{align*}
\big[ x_{ij}, x_{iij} \big] &= \big[ x_i, [x_{j}, x_{iij}] \big] 
- q_{ij} x_{j} \big[ x_i , x_{iij} \big] 
+ q_{jj}q_{ji}^2\big[ x_i , x_{iij} \big]  x_{j} = 0; \\
\big[ x_{ij}, x_{jji} \big] &= \big[ x_i, [x_{j}, x_{jji}] \big] 
- q_{ij} x_{j} \big[ x_i , x_{jji} \big] 
+ q_{ii}q_{ij}^2\big[ x_i , x_{jji} \big]  x_{j} = 0; \\
\big[ x_{iij}, x_{jji} \big] &= \big[ [x_i, x_{ij}], x_{jji} \big] \\
&= \big[ x_i, [x_{ij}, x_{jji}] \big] - q_{ii}q_{ij} x_{ij} \big[ x_i , x_{jji} \big] 
+q_{ii}^2 q_{ij}\big[ x_i , x_{jji} \big]  x_{ij} \\
& = \big[ x_i, [x_{ij}, x_{jji}] \big] = 0. 
\end{align*}

\ref{item:A-2-N=3-central-normal} We claim that the generators of $\cZ$ are annihilated by the braided adjoint action of $\htoba$. 
Fix $i \in \I_2$. By definition
$(\ad_c x_i)z_2=0=(\ad_c x_i)z_3$. In $T(V)$ we have 
$(\ad_cx_i) x_i^3= x_i^4 - q_{ii}^3 x_i^4 = 0$, and if $j\neq i$ then 
\begin{align}\label{eq:xjjji-(adxi)xj^3}
\begin{aligned}
x_{jjji} = (\ad_c x_j)^3 x_i
&= \sum_{k=0}^3 (-1)^k q_{ji}^k q_{jj}^{k(k-1)/2} \binom{3}{k}_{q_{jj}} x_j^{3-k} x_i x_j^k \\
&= x_j^3 x_i - q_{ji}^3 x_i x_j^3 = -  q_{ji}^3 (\ad_c x_i) x_j^3.
\end{aligned}
\end{align}
Thus $\ad_c x_i$ annihilates $z_1$ and $z_4$. Finally, we proceed with $z_5$.
From  \ref{item:A-2-N=3-eminent-computations} we get the commutation $x_{112}x_{12}=q_{11}^2q_{12} x_{12} x_{112}$ in $\htoba$. 
Then using \eqref{eq:braided-commutator-right-mult}
\begin{align}\label{eq:(adx1)x12^3}
\begin{aligned}
(\ad_c x_1)z_5 &
=x_{112}x_{12}^2 +q_{11}q_{12}x_{12}x_{112}x_{12}+q_{11}^2q_{12}^2 x_{12}^2 x_{112} \\
&=q_{12}^2(q_{11}^4+q_{11}^3+q_{11}^2)  x_{12}^2 x_{112} = 0.
\end{aligned}
\end{align}
For $(\ad_c x_2)z_5$, notice that on the one hand
\begin{align}\label{eq:(adx2)x12^3}
\begin{aligned}
\big[x_{12},  - \big]^3 x_2
&= \sum_{k=0}^3 (-1)^k (q_{12}q_{22})^k q_{22}^{k(k-1)/2} \binom{3}{k}_{q_{22}} x_{12}^{3-k} x_2 x_{12}^k \\
&= x_{12}^3 x_i - q_{12}^3 x_i x_{12}^3 =  - q_{12}^3 (\ad_c x_2) x_{12}^3.
\end{aligned}
\end{align}
On the other hand, using  $[x_{12}, x_{2}]= q_{12}^2q_{22} x_{221}$ 
and \ref{item:A-2-N=3-eminent-computations} 
we get
\begin{align}\label{eq:[x12,-]^3x2}
\Big[ x_{12},\big[ x_{12}, [x_{12}, x_{2}] \big] \Big] 
= q_{12}^2q_{22} \Big[ x_{12}, [x_{12}, x_{221}] \Big] = \Big[ x_{12}, 0 \Big] = 0,
\end{align}
so $(\ad_c x_2)z_5 = 0$.
This shows that $\cZ$ is a normal subalgebra.

\medbreak
Next we verify that $\Delta (z_i) \in \cZ \ot \cZ$ for $i\in \I_5$. This is clear for $i \in \I_4$, because those elements are primitive in $T(V)$; for $i=5$ we compute in $T(V)$:
\begin{align}\label{eq:A2-N=3-comult-x12^3}
\begin{aligned}
\Delta(x_{12}^3) 
=& x_{12}^3 \ot 1 + 1 \ot x_{12}^3 \\
&+(q^{-1}-q^{-2}) x_{112}\ot x_{221} 
+ (1-q^{-1})^3 q_{21}^{3} x_{1}^3\ot x_2^3\\
&+ (1-q)^2 q_{21}^{3} x_{1112}\ot x_2^2 
- (1-q^{-1})^2q^{-1} x_1^2 \ot x_{2221} \\
&+(q-1) x_1 \ot [x_{12}, x_{221}] 
- (1-q^{-1}) q_{21} [x_{12}, x_{112}] \ot x_2.
\end{aligned}
\end{align}
Using \ref{item:A-2-N=3-eminent-computations} and the defining relations of $\htoba$ we see that $\cZ$ is a Hopf subalgebra.

\medbreak
\ref{item:A-2-N=3-central-basis} We show that any pair of generators of $\cZ$ $q$-commute. 
By definition of $\htoba$, both $x_1$ and $x_2$ 
$q$-commute with $z_2$ and $z_3$, so $z_4$ and $z_1$ 
$q$-commute with $z_2$ and $z_3$. 
Secondly, \eqref{eq:xjjji-(adxi)xj^3} implies that $z_1$ and $z_4$ $q$-commute. Thirdly, \ref{item:A-2-N=3-eminent-computations} shows that $z_5$ $q$-commutes with $z_3$ and $z_2$, and also that $z_2$ and $z_3$ $q$-commute. Lastly, $z_5$ $q$-commutes with $z_4$ by \eqref{eq:(adx1)x12^3}, and with $z_1$ by \eqref{eq:(adx2)x12^3} and \eqref{eq:[x12,-]^3x2}.
Hence $B$ linearly generates $\cZ$. 

The linear independence is proven by steps.

\begin{stepooo} \label{step:A2-N=3-central-basis-1}
The set $ \left\{ z_1^{n_1} z_2^{n_2} z_3^{n_3} z_4^{n_4}  \colon n_i \in \N_0 \right\}$ is linearly independent.
\end{stepooo}
\pf
Consider the Hopf algebra $\htoba \# \ku \Z^2$; let $A$ denote the subalgebra generated by $z_1, \dots, z_4$ and $\Z^2$. 
Since all the generators of $A$ are either skew-primitives or group-likes, it follows that $A$ itself is a pointed Hopf algebra. 
Notice that $z_1, \dots, z_4 \in \mP(\cZ)$ are linearly independent. Indeed, they are non-zero because their $\Z$-degree is $3$, 
so they are linearly independent since their $\Z^2$-degrees are pairwise different 
(here we are using that the defining ideal of $\htoba$ is a Hopf ideal generated by $\Z^2$-homogeneous elements of $\Z$-degree $4$). 
Hence the infinitesimal braiding of $A$ contains the braided vector space $\ku z_1 \oplus \dots\oplus \ku z_4$, which is quantum linear space with all points labeled by $1$. Thus $ \left\{ z_1^{n_1} z_2^{n_2} z_3^{n_3} z_4^{n_4} g \colon n_i \in \N_0, g \in \Z^2 \right\} \subset A$ is linearly independent.
\epf

\begin{stepooo}\label{step:A2-N=3-central-basis-2}
The element $z_5$ does not belong to the left ideal $\htoba\langle z_1, z_2, z_3, z_4\rangle$.
\end{stepooo}
\pf
We verify this using \cite{GAP}.
\epf

The ideal $\htoba\langle z_1, z_2, z_3, z_4\rangle$ is a Hopf ideal because the generators are primitive. Denote the quotient by $R$ and consider the projection $\pi_R \colon \htoba \twoheadrightarrow R$.

\begin{stepooo}\label{step:A2-N=3-central-basis-3}
The set $ \left\{ \pi_R(z_5)^{n}  \colon n \in \N_0 \right\}$ is linearly independent.
\end{stepooo}
\pf
Consider the Hopf algebra $R \# \Z^2$. The subalgebra 
generated by $\pi_R(z_5)$ and $\Z^2$ is a pointed Hopf algebra. 
Moreover, its infinitesimal braiding contains $\pi_R(z_5)$, 
which is a non-zero point by Step \ref{step:A2-N=3-central-basis-2} and is labeled by $1$. 
Now proceed as in the proof of Step \ref{step:A2-N=3-central-basis-1}.
\epf

\begin{stepooo}\label{step:A2-N=3-central-basis-4}
We have $(\id \ot \pi_R) \Delta (z_5^n) = \sum_{k=0}^n \binom{n}{k} z_5^k \ot \pi_R(z_5)^{n-k}$ 
for all $n\in\N_0$.
\end{stepooo}
\pf
The case $n=0$ is obvious, and $n=1$ follows from \eqref{eq:A2-N=3-comult-x12^3}. An standard inductive argument for braided comultiplication yields the desired result.
\epf

\begin{stepooo}\label{step:A2-N=3-central-basis-5}
The set $B$ is linearly independent. 
\end{stepooo}
\pf
Let $\sum_{n_1, \dots,n_5 \in \N_0} \lambda_{n_1, \dots,n_5} z_1^{n_1} z_2^{n_2} z_3^{n_3} z_4^{n_4} z_5^{n_5} = 0$. Assume there exists $n_5$ such that $\lambda_{n_1, \dots,n_5}\ne 0$ for some $n_1, \dots, n_4 \in \N_0$; take $N$ as the maximal one. 
By Step \ref{step:A2-N=3-central-basis-3} there is a linear map $f\colon R \to \ku$ such that $f(\pi_R(z_5)^n) = \delta_{n,N}$ for all $n\in\N_0$. Now using Step \ref{step:A2-N=3-central-basis-4} we compute
\begin{align*}
0
&=(\id \ot f) (\id \ot \pi_R) \Delta \left(\sum_{n_1, \dots,n_5 \in \N_0}
\lambda_{n_1, \dots,n_5} 
z_1^{n_1} z_2^{n_2} z_3^{n_3} z_4^{n_4} z_5^{n_5} \right)
\\
&=\sum_{n_1, \dots,n_5 \in \N_0} \lambda_{n_1, \dots,n_5} 
(\id \ot f) (\id \ot \pi_R) 
\left( \left( \prod_{i=1}^{4} \sum_{j=0}^{n_i} \binom{n_i}{j}z_i^{j} \ot z_i^{n_i-j} \right)\Delta z_5^{n_5} \right) 
\\
&=\sum_{n_1, \dots,n_4 \in \N_0, \ n_5\leq N} \lambda_{n_1, \dots,n_5} 
(\id \ot f)
\left( \sum_{j=0}^{n_5} z_1^{n_1} z_2^{n_2} z_3^{n_3} z_4^{n_4} z_5^{j} \ot  \pi_R(z_5^{n_5-j}) \right)
\\
&=\sum_{n_1, \dots,n_4 \in \N_0} \lambda_{n_1, \dots,n_4, N} z_1^{n_1} z_2^{n_2} z_3^{n_3} z_4^{n_4} \ot 1.
\end{align*}
This contradicts Step \ref{step:A2-N=3-central-basis-1}.
\epf

\medbreak
\ref{item:A-2-N=3-central-coinvariants} Since $\Delta(z_i) \in \cZ\ot \cZ$ and $\cZ$ is normal, the right ideal $\htoba \cZ ^+$ is a Hopf ideal. By \cite[Proposition 3.6 (c)]{A+} we get that $\cZ = {}^{\co \pi }\htoba$ is equivalent to $\toba_{\bq} \simeq \htoba/\htoba \cZ ^+$. This last equality holds because the diagram
\begin{align*}
\begin{aligned}
\xymatrix@C+8pt { 
\cJ_{\bq}   \ar@{^{(}->}[r] \ar @{->>}[d] & 
T(V) \ar @{->>}[d] \ar @{->>}[dr] & 
\\
\htoba \cZ ^+  \ar@{^{(}->}[r] & 
\htoba \ar @{->>}[r]^{\pi}&
\toba_{\bq}, }
\end{aligned}
\end{align*}
commutes.
\epf

\begin{pro}\label{pro:A-2-N=3-eminent} 
\begin{enumerate}[leftmargin=*,label=\rm{(\alph*)}]
\item \label{item:A-2-N=3-eminent-extension} 
There is an extension of braided Hopf algebras
\begin{align*}
\ku \to \cZ \hookrightarrow \htoba \twoheadrightarrow \toba_{\bq} \to \ku. 
\end{align*}
\item \label{item:A-2-N=3-eminent-GK} The pre-Nichols algebra $\htoba$ is eminent and $\GK \htoba = 5$.
\end{enumerate}
\end{pro}

\pf
\ref{item:A-2-N=3-eminent-extension} Follows from Lemma \ref{lem:A-2-N=3-central} \ref{item:A-2-N=3-central-coinvariants}.

\ref{item:A-2-N=3-eminent-GK} We know that $\htoba$ covers all elements of $\pref$ by Lemma \ref{lem:A2-N=3-discard}; it remains to show that $\htoba$ itself belongs to $\pref$. By \cite[Proposition 3.6 (d)]{A+} there is a right $\cZ$-linear isomorphism $\toba_{\bq} \ot \cZ \simeq \htoba$. Since $\toba_{\bq}$ is finite dimensional, this implies that $\htoba$ is finitely generated as a $\cZ$-module.  Now \cite[Proposition 5.5]{KL} provides $\GK \htoba = \GK \cZ = 5$.
\epf

\subsection{Type $B_2$} \label{subsec:Cartan-B2}

\begin{lemma}\label{lem:eminent-B2}
Assume that $\ba$ is of Cartan type $B_2$. Then the distinguished pre-Nichols algebra $\wtoba_{\bq}$ is eminent.
\end{lemma}

\pf Here $N>2$. We may fix a braiding matrix $\bq$ such that $q_{11}=q_{22}^2$, so $q= q_{22}$ and $\widetilde{q_{12}}=q^{-2}$. 
Let $\toba$ be a finite GK-dimensional pre-Nichols algebra of $V$.
It is enough to prove that $x_{112}=0=x_{2221}$ in $\toba$. 

\smallbreak
Assume first $x_{112}\neq0$, and consider the 3-dimensional braided subspace $W:=\ku x_1 \oplus \ku x_2\oplus \ku x_{112} \subset \mP(\toba)$. Then $\GK \toba(W)<\infty$ from Lemma \ref{lem:subspace-of-primitives}. We split the proof according to the several possibilities for $N$.

\smallbreak
\noindent $\heartsuit$ $N=3$. Now the braiding on $W$ is of Cartan type 
\begin{align*}
\xymatrix{ \underset{1}{\overset{q^2}{\circ}} \ar  @{-}[rr]^{q^{-2}}  & & \underset{2}{\overset{q}{\circ} } \ar  @{-}[rr]^{q^{-2}}  & & \underset{112}{\overset{q^2}{\circ} }},
&&A= \begin{pmatrix} 2 & -1 & 0\\ -2&2&-2\\0&-1 & 2 \end{pmatrix}.
\end{align*}
Since $A$ is of affine type $C_2^{(1)}$, this contradicts \cite[Theorem 1.2(a)]{AAH-diag}.

\smallbreak
\noindent $\heartsuit$ $N=6$. In this case the braiding on $W_2:=\ku x_2\oplus \ku x_{112} \subset W$ is
\begin{align*}
\xymatrix{ \underset{ 2 }{\overset{q}{\circ}} \ar  @{-}[rr]^{q^{-2}}  & & \underset{ 112 }{\overset{q^5}{\circ} }}, && \text{Cartan type }
\begin{pmatrix} 2 & -2 \\ -4 & 2 \end{pmatrix}.
\end{align*}
The Cartan matrix is of indefinite type, and by \cite[Theorem 1.2(b)]{AAH-diag} this contradicts $\GK \toba(W)<\infty$ .

\smallbreak
\noindent $\heartsuit$ $N\neq 3, 6$. The Dynkin diagram of $W_1:=\ku x_1\oplus \ku x_{112} \subset W$ is  
$$\bD_1= \xymatrix{ \underset{1}{\overset{q^2}{\circ}} \ar  @{-}[rr]^{q^6}  & & \underset{112}{\overset{q^5 }{\circ}}}.$$
Since $\GK \toba(W_1)<\infty$, it follows from \cite[Theorem 1.2(b)]{AAH-diag} that the associated root system is finite. Now $\bD_1$ is connected; by exhaustion on \cite[Table 1]{H-class}, we deduce that it must be $N=4$ or $N=8$. We turn again to $W_2$, whose Dynkin diagram is easily computed in each case:
\begin{align*}
\heartsuit\heartsuit N=4: &&\xymatrix{ \underset{ 2 }{\overset{q}{\circ}} \ar  @{-}[rr]^{q^{-2}}  & & \underset{ 112 }{\overset{q}{\circ} }}, && \text{Cartan type }
\begin{pmatrix} 2 & -2 \\ -2 & 2 \end{pmatrix},
\\
\heartsuit\heartsuit N=8: && \xymatrix{ \underset{ 2 }{\overset{q}{\circ}} \ar  @{-}[rr]^{q^{-2}}  & & \underset{ 112 }{\overset{q^5 }{\circ} }}, &&
\text{Cartan type }
\begin{pmatrix} 2 & -2\\ -2 & 2 \end{pmatrix}.
\end{align*}
In any case the Cartan matrix is of affine type $A_1^{(1)}$, so $\GK \toba(W_2) = \infty$ by \cite[Theorem 1.2(b)]{AAH-diag}.

\medbreak
Assume $x_{2221}\neq0$ in $\toba$. The subspace $U :=\ku x_1 \oplus \ku x_2\oplus \ku x_{2221} \subset \mP(\toba)$ has dimension 3 and $\GK \toba(U)<\infty$. Now $U_1:=\ku x_1\oplus \ku x_{2221} \subset U$ has connected Dynkin diagram  
$$\xymatrix{ \underset{1}{\overset{q^2}{\circ}} \ar  @{-}[rr]^{q^2}  & & \underset{2221}{\overset{q^5 }{\circ}}},$$
and it is finite by \cite[Theorem 1.2(b)]{AAH-diag}. By exhaustion on \cite[Table 1]{H-class} we deduce that $N=4$. Then the Dynkin diagram of $U$ is of Cartan type 
\begin{align*}
\xymatrix{ \underset{2}{\overset{q}{\circ}} \ar  @{-}[rr]^{-1}  & & \underset{1}{\overset{-1}{\circ} } \ar  @{-}[rr]^{-1}  & & \underset{2221}{\overset{q}{\circ} }},
&&A= \begin{pmatrix} 2 & -2 & 0\\ -1&2&-1\\0&-2 & 2 \end{pmatrix}.
\end{align*}
Since $A$ is of affine type $C_2^{(1)}$, this contradicts \cite[Theorem 1.2(a)]{AAH-diag}.
\epf

\subsection{Type $G_2$} \label{subsec:Cartan-G2}

\begin{lemma}\label{lem:preNichols-G2=image-distinguished}
Assume that $\ba$ is of Cartan type $G_2$. Then the quantum Serre relations hold in any $\toba \in \pref$. In particular, $\wtoba_{\bq}$ is eminent if $N\neq 4, 6$.
\end{lemma}

\pf
Here $N>3$. Let $\toba \in \pref(V)$; we show first that the quantum Serre relations $x_{11112}=0=x_{221}$ hold in $\toba$. 

\smallbreak
Start assuming $x_{11112}\neq0$. Then the 3-dimensional subspace $W:=\ku x_1 \oplus \ku x_2\oplus \ku x_{11112} \subset \mP(\toba)$ satisfies $\GK \toba(W)\leq \GK \toba$ by Lemma  \ref{lem:subspace-of-primitives}. 
The Dynkin diagram of $W_1:=\ku x_1\oplus \ku x_{11112} \subset W$ is
$$\bD_1= \xymatrix{ \underset{1}{\overset{q}{\circ}} \ar  @{-}[rr]^{q^5}  & & \underset{11112}{\overset{q^7 }{\circ}}}.$$
Since $\GK \toba(W_1)<\infty$, it follows from \cite[Theorem 1.2(b)]{AAH-diag} that the root system of $\bD_1$ is finite. We split the proof according to the several possibilities for $N$.

\smallbreak
\noindent $\heartsuit$ $N=5$. The diagram $\bD_1$ is disconnected, but we might consider instead $W_2:= \ku x_{11112}\oplus \ku x_2\subset W$, that satisfies $\GK \toba(W_2)<\infty$ as well. By direct computation $W_2$ is of indefinite Cartan type:
\begin{align*}
\xymatrix{ \underset{11112}{\overset{q^2}{\circ}} \ar  @{-}[rr]^{q^{-6}}  & & \underset{2}{\overset{q^3}{\circ} }}, &&
\begin{pmatrix} 2 & -3 \\ -2 & 2 \end{pmatrix},
\end{align*}
which is in contradiction with \cite[Theorem 1.2(a)]{AAH-diag}.

\smallbreak
\noindent $\heartsuit$ $N\neq 5$. Now $\bD_1$ is connected and finite; by inspection on \cite[Table 1]{H-class}, it must be $N=4$ or $N=6$.

$\heartsuit\heartsuit N=4$.
In this case $W_2$ is of Cartan type
\begin{align*}
\xymatrix{ \underset{11112}{\overset{q^3}{\circ}} \ar  @{-}[rr]^{q^{-6}}  & & \underset{2}{\overset{q^3}{\circ} }}, &&
\begin{pmatrix} 2 & -2 \\ -2 & 2 \end{pmatrix},
\end{align*}
which is of affine type $A_1^{(1)}$, now contradicting \cite[Theorem 1.2(b)]{AAH-diag}.

$\heartsuit\heartsuit N=6$. 
In this case the Dynkin diagram of $W$ is of Cartan type 
\begin{align*}
\xymatrix{ \underset{11112}{\overset{q}{\circ}} \ar  @{-}[rr]^{q^{-1}}  & & \underset{1}{\overset{q}{\circ} } \ar  @{-}[rr]^{q^{-3}}  & & \underset{2}{\overset{q^3}{\circ} }},
&&A=\begin{pmatrix} 2 & -1 & 0\\ -1&2&-3\\0&-1 & 2 \end{pmatrix}.
\end{align*}
By \cite[Theorem 1.2(b)]{AAH-diag} this contradicts $\GK \toba(W)  < \infty$, since $A$ is of affine type $G_2^{(1)}$.

\medspace
Assume now $x_{221}\neq0$ in $\toba$. The subspace $U :=\ku x_1 \oplus \ku x_2\oplus \ku x_{221} \subset \mP(\toba)$ has dimension 3 and $\GK \toba(U)<\infty$. Consider two possibilities for $N$.

\smallbreak
\noindent
$\heartsuit N\neq4$. Now $U_1:=\ku x_1\oplus \ku x_{221} \subset U$ has connected Dynkin diagram  
$$\xymatrix{ \underset{1}{\overset{q}{\circ}} \ar  @{-}[rr]^{q^{-4}}  & & \underset{221}{\overset{q^7 }{\circ}}}.$$
By exhaustion on \cite[Table 1]{H-class} we conclude that this diagram is never finite, which contradicts \cite[Theorem 1.2(b)]{AAH-diag}, as $\GK \toba(U_1)<\infty$.

\smallbreak
\noindent
$\heartsuit N=4$. In this case the braiding on $U$ is of Cartan type 
\begin{align*}
\xymatrix{ \underset{1}{\overset{q}{\circ}} \ar  @{-}[rr]^{q^{-3}}  & & \underset{2}{\overset{q^{3}}{\circ} } \ar  @{-}[rr]^{q^{-3}}  & & \underset{221}{\overset{q^{3}}{\circ} }},
&&A= \begin{pmatrix} 2 & -3 & 0\\ -1&2&-1\\0&-1 & 2 \end{pmatrix}.
\end{align*}
Since $A$ is of affine type  $G^{(1)}_2$, this contradicts \cite[Theorem 1.2(a)]{AAH-diag}.

\smallbreak
Thus the quantum Serre relations hold in $\toba$. By Remark \ref{rem:qsr-presentation} this proves the assertion regarding $N\neq 4, 6$.
\epf

\subsection{Type $A_3$} \label{subsec:typeA3-eminent}

\begin{lemma}\label{lem:eminent-A3-N>2}
If $\ba$ is of Cartan type $A_3$ with $N>2$, then $\wtoba_{\bq}$ is eminent.
\end{lemma}
\pf
As $N>2$, the ideal $\cI_{\bq}$ is generated by the quantum Serre relations
$x_{13}=0$ and $x_{iij}=0$ for $|j-i|=1$, cf. \cite[p. 397]{AA-diag}.
Let $\toba \in \pref(\bq)$. Then $x_{13}=0$ holds in $\toba$ since the braided vector space $\ku x_1\oplus \ku x_3$ satisfies the hypothesis in Proposition \ref{prop:xij=0-preNichols}. 

Turn to $x_{iij}$ for some fix $i,j \in \I_3$ with $|j-i|=1$; in this case $\ku x_i\oplus \ku x_j$ is of Cartan type $A_2$. If $N>3$, then $x_{iij}=0$ in $\toba$ by Lemma \ref{lem:preNichols-A2-N>3-eminent}. It only remains the case $N=3$. Now we have $\bq(2\alpha_i+\alpha_j, 2\alpha_i+\alpha_j )=q_{ii}^4\widetilde{q_{ij}}^2q_{jj}=q^5q^{-2}=1$. Using \cite[Lemma 2.8] {AAH-triang}, in order to guarantee $x_{iij}=0$ in $\toba$ it is enough to find $k \in \I_3$ such that $\wbq (\alpha_k, 2\alpha_i+\alpha_j) \neq 1$. It is straightforward to verify that the unique $k \in \I_3$ different from $i$ and $j$ does the trick.
\epf
\subsection{Types $B_3$ and $C_3$} \label{subsec:typesB3-C3-eminent}
\begin{lemma}\label{lem:eminent-B3-C3} The distinguished pre-Nichols algebra $\wtoba_{\bq}$ is eminent
if  either 
\begin{multicols}{2}
\begin{enumerate}[leftmargin=*,label=\rm{(\roman*)}]
\item\label{item:eminent-B3} $\ba$ is of type $B_3$, or
\item\label{item:eminent-C3} $\ba$ is of type $C_3$.
\end{enumerate}
\end{multicols}
\end{lemma}
\pf
Let $\toba \in \pref(\bq)$. Then $x_{13}=0$ holds in $\toba$. Indeed, the braided vector space $\ku x_1\oplus \ku x_3$ satisfies the hypothesis in Proposition \ref{prop:xij=0-preNichols}. Similarly, since $\ku x_2 \oplus\ku x_3$ is of type $B_2$, it follows from Lemma \ref{lem:eminent-B2} that the quantum Serre relations involving $x_2$ and $x_3$ hold in $\toba$.

\begin{stepo}\label{step:B3-eminent-qsr-hold}
If $\ba$ is of Cartan type $B_3$, then the quantum Serre relations hold in any finite $\GK$ pre-Nichols algebra.
\end{stepo}
\pf
Here $\ku x_1 \oplus\ku x_2$ has Dynkin diagram $\xymatrix{\overset{q^2}{\circ}\ar  @{-}[r]^{q^{-2}}  &\overset{q^2}{\circ }}$, type $A_2$. Hence, if $\ord q^2 > 3$, we know from Lemma \ref{lem:preNichols-A2-N>3-eminent} that the quantum Serre relations between $x_1$ and $x_2$ hold in $\toba$. Let us show that in the cases $\ord q^2 =2, 3$ the same happens.

\smallbreak
\noindent $\heartsuit$ $\ord q^2 =2$. If $x_{112} \neq0$ in $\toba$, we get a subspace $\ku x_2 \oplus \ku x_3 \oplus \ku x_{112} \subset \mP(\toba)$ of dimension $3$ with the following Dynkin diagram
\begin{align*}
\xymatrix{ \underset{2}{\overset{-1}{\circ}} \ar  @{-}[rr]^{-1}  & & \underset{3}{\overset{q}{\circ} } \ar  @{-}[rr]^{-1}  & & \underset{112}{\overset{-1}{\circ}} },
&& \text {Cartan type } \begin{pmatrix} 2 & -1 & 0\\ -2&2&-2\\0&-1 & 2 \end{pmatrix}.
\end{align*}
This matrix is of affine type $C_2^{(1)}$, hence $\GK \toba = \infty$, a contradiction.
\\
Similarly, the assumption $x_{221} \neq0$ yields a subspace of $\mP(\toba)$ with braiding
\begin{align*}
\begin{aligned}
\xymatrix@C+24pt { 
& \overset{-1}{\underset{221}{\circ}} \\
\overset{-1}{\underset{1}{\circ}} \ar  @{-}[r]^{-1}  & \overset{-1}{\underset{2}{\circ}} \ar  @{-}[u]^{-1}_{}  \ar  @{-}[r]^{-1} &
\overset{q}{\underset{3}{\circ}},}
\end{aligned}
& &\text {Cartan type } 
\begin{pmatrix} 2&-1&0&0\\ -1&2&-1&-1 \\0&-2&2&0 \\0&-1&0&2  \end{pmatrix}. 
\end{align*}
The Cartan matrix is of affine type $B_3^{(1)}$, and again $\GK (\toba) = \infty$.
\smallbreak
\noindent $\heartsuit$ $\ord q^2 =3$. Notice that 
\begin{align*}
\bq(2\alpha_1+\alpha_2, 2\alpha_1+\alpha_2)&=q^6=1, & 
\wbq(2\alpha_1+\alpha_2, \alpha_3)&=q^{-2}\neq 1;
\\ 
\bq(\alpha_1+2\alpha_2, \alpha_1+2\alpha_2)&=q^6=1, & 
\wbq(\alpha_1+2\alpha_2, \alpha_3)&=q^{-4}\neq 1. 
\end{align*} 
Assuming $x_{112}\neq0$ in $\toba$ we get $\ku x_3 \oplus \ku x_{112} \subset \mP(\toba)$ with Dynkin diagram $\xymatrix{\overset{q}{\circ}\ar  @{-}[r]^{q^{-2}}  &\overset{1}{\circ }}$. Then by \cite[Lemma 2.8]{AAH-triang} it follows $\GK \toba =\infty$, a contradiction. By  the same argument, it can not be $x_{221}\neq 0$ in $\toba$.
\epf

The assertion \ref{item:eminent-B3}  for $N>4$ follows since, in that case,  $\wtoba_{\bq}$ is presented by the quantum Serre relations, cf. Remark \ref{rem:qsr-presentation}.

\begin{stepo}\label{step:B3-N=3-eminent}
If $\ba$ is of Cartan type $B_3$ with $N=3$, then $\wtoba_{\bq}$ is eminent.
\end{stepo}
\pf
By \cite[pp. 399, 400]{AA-diag}, $\wtoba_{\bq}$ is presented by the quantum Serre relations and $[x_{3321}, x_{32}]_c=0$. Given $\toba \in \pref$, let us show that $[x_{3321}, x_{32}]_c \in \mP(\toba)$. Using  $x_{13}=0$ an straightforward computation gives
\begin{align*}
\Delta(x_{3321}) =& x_{3321}\ot 1 + 1\ot x_{3321} + (1-q_{33})x_{332}\ot x_1
\\
&+ (1-q_{33})q_{33} x_3 \ot x_{321} + (1-q_{33}) (1-q_{22}) x_3^2 \ot x_{21}.
\end{align*}
with this we compute 
\begin{align*}
\Delta([x_{3321}, x_{32}]_c) =& [x_{3321}, x_{32}]_c\ot 1 + 1\ot [x_{3321}, x_{32}]_c \\
&- (1-q_{33})^2(1-q_{22})q_{12} \, x_3^2\ot x_{221}\\
&- (1-q_{33})q_{12}q_{13}q_{23} q_{33} \, x_{3332} \ot x_{21}\\
&+ (1-q_{33})q_{33}^2q_{12}q_{13} \, x_{332} \ot (x_{321}-[x_{32},x_1]_c)\\
&- (1-q_{33})q_{23}q_{13} \, x_{33321} \ot x_2\\
&+ (1-q_{33})q_{13}q_{12}\, [x_{332},x_{32}]_c \ot x_1 \\
&+ (1-q_{33})q_{33} \, x_3 \ot (q_{13}q_{23}q_{33} [x_{3321},x_2]_c +[x_{321},x_{32}]_c)
\end{align*}
The third and fourth terms vanishes in $\toba$ by Step \ref{step:B3-eminent-qsr-hold}. For the fifth term, an straightforward computation involving $x_{13}=0$ shows that $x_{321}=[x_{32},x_1]_c$. The last three terms also vanish, but they require a more detailed analysis.

\smallbreak
\noindent $\heartsuit$ $x_{33321}=0$ in $\toba$. Notice that
\begin{align*}
\Delta(x_{33321}) =& x_{33321}\ot 1 + 1\ot x_{33321} \\
&+ (1-q_{33})x_{3332}\ot x_1 - (1-q_{33}^2)q_{32} x_{332} \ot x_{31},
\end{align*}  
so this element is primitive in $\toba$ by Step \ref{step:B3-eminent-qsr-hold}. Assuming $x_{33321}\neq0$ we get a subspace $\ku x_1 \oplus \ku x_{33321} \subset \mP(\toba)$ where the braiding is given by 
\begin{align*}
\xymatrix{ \underset{1}{\overset{q^2}{\circ}} \ar  @{-}[rr]^{q^{-1}}  & & \underset{33321}{\overset{q^2}{\circ} }}, && \text{ Cartan type }
\begin{pmatrix} 2 & -2 \\ -2 & 2 \end{pmatrix}, && \text{ affine type } A_1^{(1)}.
\end{align*}
this contradicts $\GK \toba <\infty$ by \cite[Theorem 1.2]{AAH-diag}.

\smallbreak
\noindent $\heartsuit$ $[x_{332},x_{32}]_c=0$ in $\toba$. Now we have 
\begin{align*}
\Delta([x_{332}, x_{32}]_c) =& [x_{332}, x_{32}]_c\ot 1 + 1\ot [x_{332}, x_{32}]_c \\
&- (1-q_{33})^2q_{23} \, x_3^2\ot [x_{32}, x_{2}]_c - (1-q_{33})q_{33}q_{23} \, x_{3332} \ot x_2.
\end{align*}
The element $[x_{32}, x_2]_c$ is primitive in $\toba$, so it vanishes by the same reason that $x_{223}$ does (cf. proof of Lemma \ref{lem:eminent-B2}). So $[x_{332}, x_{32}]_c \in \mP(\toba)$ by Step \ref{step:B3-eminent-qsr-hold}. If it is non-zero, consider $\ku x_1 \oplus \ku[x_{332},x_{32}]_c \subset \mP(\toba)$ where the braiding is 
\begin{align*}
\xymatrix{ \underset{x_1}{\overset{q^2}{\circ}} \ar  @{-}[rr]^{q^{-1}}  & & \underset{[x_{332}, x_{32}]}{\overset{q^2}{\circ} }}, && \text{ Cartan type }
\begin{pmatrix} 2 & -2 \\ -2 & 2 \end{pmatrix}, && \text{ affine type } A_1^{(1)},
\end{align*}
thus we get the same contradiction as with $x_{33321}$.

\smallbreak
\noindent $\heartsuit$ $q_{13}q_{23}q_{33} [x_{3321},x_2]_c +[x_{321},x_{32}]_c=0$. Denote this element by $r$. Then
\begin{align*}
\Delta(r) =& r \ot 1 + 1\ot r + (1-q_{33})q_{22}q_{12}q_{13} \, x_{32}\ot (x_{321}-[x_{32}, x_{1}])\\
& + (1-q_{33})q_{33}q_{12}q_{13}q_{23} \, \big[x_3, [x_{32},x_2]\big] \ot x_1\\
&- (1-q_{33})q_{22}q_{12}q_{13}q_{23} \, x_3 \ot \big[ [x_{32},x_2], x_1\big].
\end{align*}
Since $[x_{32},x_2]=0$ and $x_{321}-[x_{32}, x_{1}]=0$ in $\toba$, it follows that $r$ is primitive. If $r\neq0$ we consider $\ku x_2 \oplus \ku r \subset \mP(\toba)$. The Dynkin diagram is computed: 
\begin{align*}
\xymatrix{ \underset{2}{\overset{q^2}{\circ}} \ar  @{-}[rr]^{q^{-1}}  & & \underset{r}{\overset{q^2}{\circ} }}, && \text{ Cartan type }
\begin{pmatrix} 2 & -2 \\ -2 & 2 \end{pmatrix}, && \text{ affine type } A_1^{(1)},
\end{align*}
thus we get the same contradiction as before.

\medbreak
Using this three $\heartsuit$ we get $[x_{3321}, x_{32}]_c \in \mP(\toba)$. If this element is non-zero, consider $U=\ku x_3 \oplus \ku [x_{3321}, x_{32}]_c\subset \mP(\toba)$. We compute the braiding: 
\begin{align*}
\bq(\alpha_1+2\alpha_2+3\alpha_3,\alpha_1+2\alpha_2+3\alpha_3) &=1, && \wbq(\alpha_1+2\alpha_2+3\alpha_3,\alpha_3) = q^{-1}\neq 1.
\end{align*}
From \cite[Lemma 2.8] {AAH-triang} it follows $\GK \toba(U)=\infty$, but this contradicts $\GK\toba < \infty$. Then $[x_{3321}, x_{32}]_c=0$ in $\toba$ and Step \ref{step:B3-N=3-eminent} holds.
\epf
\begin{stepo}\label{step:B3-N=4-eminent}
If $\ba$ is of Cartan type $B_3$ with $N=4$, then $\wtoba_{\bq}$ is eminent.
\end{stepo}
\pf
By \cite[pp. 399, 400]{AA-diag}, $\wtoba_{\bq}$ is presented by the quantum Serre relations and $[x_{123}, x_2]_c=0$. We claim that this element is primitive in any $\toba \in \pref$. Indeed, using that $x_{13}=0$ in $\toba$, we get 
\begin{align*}
\Delta([x_{123}, x_2]_c) =& [x_{123}, x_2]_c\ot 1 + 1\ot [x_{123}, x_2]_c 
\\
&-(1-\widetilde{q_{12}})q_{32} x_1 \ot x_{223} + (1-\widetilde{q_{23}})q_{32} [x_{12}, x_2]_c \ot x_3.
\end{align*}
By straightforward computations, $[x_{12}, x_2]_c \in \mP (\toba)$ and it vanishes by the same reason that $x_{221}$ does (cf. proof of Lemma \ref{lem:preNichols-A2-N>3-eminent}). Since $x_{223}=0$, the claim follows.

Assume $[x_{123}, x_2]_c\neq0$. Inside $\mP (\toba)$ we have the 2-dimensional subspace $U=\ku x_3 \oplus \ku [x_{123}, x_2]_c$ where the braiding is given by
\begin{align*}
\xymatrix{ \underset{x_3}{\overset{q}{\circ}} \ar  @{-}[rr]^{-1}  & & \underset{[x_{123}, x_2]}{\overset{-q}{\circ} }}, && \text { Cartan type }
\begin{pmatrix} 2 & -2 \\ -2 & 2 \end{pmatrix}.
\end{align*}
Since this matrix is of affine type $A_1^{(1)}$, from \cite[Theorem 1.2(b)]{AAH-diag} it follows $\GK\toba(U)=\infty$, contradicting $\toba \in \pref$.
\epf

\begin{stepo}\label{step:C3-qsr-hold}
If $\ba$ is of Cartan type $C_3$, then the quantum Serre relations hold in any finite $\GK$ pre-Nichols algebra.
\end{stepo}
\pf
Now $\ku x_1 \oplus\ku x_2$ has Dynkin diagram $\xymatrix{\overset{q}{\circ}\ar  @{-}[r]^{q^{-1}}  &\overset{q}{\circ }}$, type $A_2$. If $N>3$, then the quantum Serre relations in $x_1$ and $x_2$ hold by Lemma \ref{lem:preNichols-A2-N>3-eminent}. For the case $N=3$, let $i, j$ such that $\{i, j\}=\{1,2\}$ and suppose $x_{iij}\neq0$ in $\toba$. Since 
$\bq(2\alpha_i + \alpha_j, 2\alpha_i + \alpha_j) = q^5 q^{-2}=1$ and $\wbq(2\alpha_i + \alpha_j, \alpha_3)=\widetilde{q_{i3}}^2\widetilde{q_{j3}}\neq 1$, we get $\GK \toba = \infty$ by \cite[Lemma 2.8] {AAH-triang}. 
\epf

The assertion \ref{item:eminent-C3}  for $N>4$ follows since, in that case,  $\wtoba_{\bq}$ is presented by the quantum Serre relations, 
see Remark \ref{rem:qsr-presentation}.

\begin{stepo}\label{step:C3-N=3-eminent}
If $\ba$ is of Cartan type $C_3$ with $N=3$, then $\wtoba_{\bq}$ is eminent.
\end{stepo}
\pf
Following \cite[pp. 401, 402]{AA-diag}) we see that $\wtoba_{\bq}$ is presented by the quantum Serre relations and $\big[ [x_{123}, x_2]_c, x_2\big]_c=0$. Given $\toba \in \pref$, let us show that this element is primitive in $\toba$. Using $x_{13}=0$ it follows
\begin{align*}
\Delta([x_{123}, x_2]_c) =& [x_{123}, x_2]_c\ot 1 + 1\ot [x_{123}, x_2]_c + (1-\widetilde{q_{12}}) x_{123} \ot x_2
\\
&+ (1 - q^2) x_{12}\ot x_{32} 
+ (1-\widetilde{q_{12}})x_1 \ot (x_{23}x_2-q_{32}x_2x_{23}) 
\\
&+ (1-\widetilde{q_{23}})q_{32} [x_{12}, x_2]_c \ot x_3.
\end{align*}
By straightforward computations, $[x_{12}, x_2]_c = q_{12}^2qx_{221}$ in $T(V)$, and so $[x_{12}, x_2]_c$ vanishes in $\toba$ by Step \ref{step:C3-qsr-hold}. Then we obtain
\begin{align*}
\Delta &\big( \big[ [x_{123}, x_2]_c, x_2\big]_c\big) = \big[ [x_{123}, x_2]_c, x_2\big]_c \ot 1 + 1\ot \big[ [x_{123}, x_2]_c, x_2\big]_c 
\\
&+ (1-q^2)q_{32}q_{22}[x_{12}, x_2]_c \ot x_{32} + (1-\widetilde{q_{12}})q_{22}q_{32}^2 x_1 \ot x_{2223},
\end{align*}
and now the claim follows from Step \ref{step:C3-qsr-hold}.

If $\big[ [x_{123}, x_2]_c, x_2\big]_c\neq0$, consider $U=\ku x_1\oplus \ku \big[ [x_{123}, x_2]_c, x_2\big]_c\subset \mP(\toba)$. 
By \cite[Lemma 2.8] {AAH-triang}, since 
$\bq(\alpha_1+3\alpha_2+\alpha_3,\alpha_1+3\alpha_2+\alpha_3)= q^{12}q^{-9} = 1$
and
$\wbq(\alpha_1+3\alpha_2+\alpha_3,\alpha_1)=q^2q^{-3}\neq 1$, we have $\GK \toba(U)=\infty$. This contradicts $\toba \in \pref$.
\epf

\begin{stepo}\label{step:C3-N=4-eminent}
If $\ba$ is of Cartan type $C_3$ with $N=4$, then $\wtoba_{\bq}$ is eminent.
\end{stepo}
\pf
By \cite[Theorem 3.1]{An-diagonal}, $\wtoba_{\bq}$ is presented by the quantum Serre relations and $[x_{123}, x_{23}]_c=0$. Let us show that this element is primitive in any finite dimensional pre-Nichols algebra $\toba$. 

First we claim that $[x_{123}, x_3]_c=0$ in $\toba$: using that $x_{13}=0$ we compute
\begin{align*}
\Delta \big( [x_{123}, x_3]_c\big) =& [x_{123}, x_3]_c \ot 1 + 1\ot [x_{123}, x_3]_c 
\\
&+ (1-\widetilde{q_{12}}q_{23}q_{33})x_{13}\ot x_{23} + (1-\widetilde{q_{12}}) x_1 \ot [x_{23},x_3]_c.
\end{align*}
Since $[x_{23}, x_3]_c \in \mP (\toba)$, it vanishes in $\toba$ by the same reason that $x_{332}$ does (cf. proof of Lemma \ref{lem:eminent-B2}). So $[x_{123}, x_3]_c \in \mP(\toba)$. Hence, if it is non-zero we get a subspace $U=\ku x_1 \oplus [x_{123}, x_3]_c \subset \mP(\toba)$ where the braiding is given by
\begin{align*}
\xymatrix{ \underset{x_3}{\overset{q}{\circ}} \ar  @{-}[rr]^{q^{-3}}  & & \underset{[x_{123}, x_3]}{\overset{q}{\circ} }}, &&\text{ indefinite Cartan type }
\begin{pmatrix} 2 & -3 \\ -3 & 2 \end{pmatrix}.
\end{align*}
But then $\GK\toba(U)=\infty$ by \cite[Theorem 1.2(b)]{AAH-diag}, a contradiction.

Next we compute
\begin{align*}
\Delta &\big( [x_{123}, x_{23}]_c\big) =[x_{123}, x_{23}]_c \ot 1 + 1\ot [x_{123}, x_{23}]_c 
\\
&+ (1-\widetilde{q_{23}})q_{12}q_{22}q_{32}x_{2}\ot [x_{123}, x_3] + (1-\widetilde{q_{23}})^2q_{32} x_1 \ot [x_{12},x_2]_c.
\end{align*}
Using the previous claim and the fact $[x_{12},x_2]_c=q_{12}^2qx_{221}=0$ (by Step \ref{step:C3-qsr-hold}), we get $[x_{123}, x_{23}]_c \in \mP(\toba)$. If $[x_{123}, x_{23}]_c \neq0$, consider the subspace $W = \ku x_1 \oplus \ku x_2 \oplus \ku [x_{123}, x_{23}]_c \subset \mP(\toba)$, where the braiding is 
\begin{align*}
\xymatrix{ \underset{1}{\overset{q}{\circ}} \ar  @{-}[rr]^{q^{-1}}  & & \underset{2}{\overset{q}{\circ} } \ar  @{-}[rr]^{q^{-1}}  & & \underset{[x_{123}, x_{23}]}{\overset{q^{3}}{\circ} }},
&& \text{ Cartan type }\begin{pmatrix} 2 & -1 & 0\\ -1&2&-1\\0&-3 & 2 \end{pmatrix}.
\end{align*}
Since the Cartan matrix is of affine type  $G^{(1)}_2$, it follows $\GK \toba(W) =\infty$ by \cite[Theorem 1.2(a)]{AAH-diag}. This contradicts $\toba \in \pref$.
\epf
The result follows.
\epf

\subsection{Some cases in rank $> 3$} \label{subsec:rank>3-eminent}
Here we assume that $\theta \ge 4$.

\begin{lemma}\label{lem:distinguished-simplylaced-eminent}
In any of the following cases, $\wtoba_{\bq}$ is eminent.

\begin{enumerate}[leftmargin=*,label=\rm{(\alph*)}]
\item\label{item:disntiguished-general1} $\ba$ is of Cartan type with simply laced Dynkin diagram and  $N>2$.

\item\label{item:disntiguished-general2} $\ba$ is of type $B_\theta$, $C_\theta$ ($\theta \ge 4$) or  $F_4$, and $N>4$.

\end{enumerate}
\end{lemma}

\pf By Remark \ref{rem:qsr-presentation} and the restrictions on $N$, $\wtoba_{\bq}$ is presented by the quantum Serre relations.
Let $\toba \in \pref(\bq)$. If $a_{ij} = 0$, then $x_{ij}=0$ holds in $\toba$ by Proposition \ref{prop:xij=0-preNichols}.
If $a_{ij} \neq 0$, then there is $k \in \I$ such that $\{i,j,k\}$ span a subdiagram of type $A_3$, $B_3$ or $C_3$. Then
$(\ad x_i)^{1-a_{ij}} (x_j)=0$ by Lemmas \ref{lem:eminent-A3-N>2} or \ref{lem:eminent-B3-C3}.  Thus $\wtoba_{\bq} \twoheadrightarrow \toba$.
\epf

In the next subsections we treat some remaining cases with small $N$.

\subsection{Types $B_{\theta}$, $C_{\theta}$, $F_4$, $\theta > 3$, $N = 3, 4$} 
\begin{lemma}\label{lem:eminent-F4}
If $\ba$ is of types $B_{\theta}$, $C_{\theta}$, with $\theta > 3$, or $F_4$, and $N = 3$ or $4$, then $\wtoba_{\bq}$ is eminent.
\end{lemma}
\pf We split the proof according to the type. Let $\toba \in \pref$.
 
\smallbreak
$\heartsuit$ Type $F_4$. Here $\wtoba_{\bq}$ is presented by the quantum Serre relations and
\begin{align} \label{eq:F4-distinguished-notqsr}
[x_{123}, x_{23}]_c=[x_{432} , x_3]_c= 0 \text { if } N=4;  &&
[x_{2234} , x_{23}]_c = 0 \text { if } N=3.
\end{align}
Since $N>2$ we get $x_{14}=0$ in $\toba$ from Proposition \ref{prop:xij=0-preNichols}. 
The subdiagram spanned by $\{1,2,3\}$ is of type $C_3$ thus the quantum Serre relations involving this indices hold in $\toba$ by Lemma 
\ref{lem:eminent-B3-C3} \ref{item:eminent-C3}.
Finally, $\{4,3,2\}$ span a diagram of type $B_3$ so the quantum Serre  relations involving this indices hold in $\toba$ by Lemma \ref{lem:eminent-B3-C3} \ref{item:eminent-B3}. Moreover \eqref{eq:F4-distinguished-notqsr} are defining relations of the distinguished pre-Nichols algebra of type $B_3$ or $C_3$ for the corresponding $N$, hence Lemma \ref{lem:eminent-B3-C3} implies that these also vanish in $\toba$.

\smallbreak
$\heartsuit$ Type $B_{\theta}$. Here $\wtoba_{\bq}$ is presented by the quantum Serre relations and
\begin{align} \label{eq:Btheta-distinguished-notqsr}
 [x_{(i \, i+2)}, x_{i+1}]_c, \, i<\theta-1, \text { if } N=4; &&
 [x_{\theta  \theta \theta-1  \theta-2}, x_{\theta  \theta-1}]_c, \, \text { if } N=3.
\end{align}
The relations involving the indices $\{\theta-2, \theta -1, \theta\}$ hold in $\toba$ by Lemma \ref{lem:eminent-B3-C3} \ref{item:eminent-B3}; also $x_{i \theta}=0$ for any $i<\theta-1$ by Proposition \ref{prop:xij=0-preNichols}. We are left to treat the relations involving $\{1, \dots, \theta-1\}$. If $N=3$ we only have the quantum Serre relations, which hold by Lemma \ref{lem:eminent-A3-N>2}. 
Turn to $N=4$. Now $\{1, \dots, \theta-1\}$ form a subdiagram of type $A_{\theta-1}$ at a root of order $2$. If $\theta -1 \geq4$ we apply Lemma \ref{lem:eminent-Atheta-N=2} to get all the Serre relations except for $x_{221}$ and $x_{\theta -2  \,\theta -2\, \theta -1}$. The last one holds by Lemma \ref{lem:eminent-B3-C3} \ref{item:eminent-B3} and the first one falls since 
the diagram
\begin{align*}
\xymatrix @C40pt @R-15pt
{ \overset{-1}{\underset{221 }{\circ}}&&&&\\
\overset{-1}{\underset{1 }{\circ}}\ar  @{-}[r]^{-1}  &
\overset{-1}{\underset{2 }{\circ}}\ar  @{-}[r]^{-1} \ar  @{-}[ul]_{-1}&
\overset{-1}{\underset{3 }{\circ}}\ar@{.}[r] & 
\overset{-1}{\underset{\theta -1}{\circ}} \ar  @{-}[r]^{-1}  & \overset{q}{\underset{\theta }{\circ}}
}
\end{align*}
is of indefinite Cartan type. Now $[x_{(i \, i+2)}, x_{i+1}]_c=0$ for $ i<\theta-2$ hold by Lemma \ref{lem:eminent-Atheta-N=2} \ref{item:Atheta-N=2-[x(i,i+2),xi+1]=0}. We treat separately the last case standing. 

$\heartsuit \heartsuit$ Type $B_4$ with $N=4$. The relations $x_{221}$ and $x_{\theta -2  \,\theta -2\, \theta -1}$ hold by the same reason as above. Moreover, we also have $x_{13}=0$. This follows from \cite[Lemma 2.8]{AAH-triang} since $\bq(\alpha_1+\alpha_3,\alpha_1+\alpha_3)=1$ and $\wbq(\alpha_1+\alpha_3,\alpha_4)\neq1$. Finally, using \eqref{eq:A3-N=2-comult-[x(13),x2]} and the relations deduced so far, we get that $[x_{(13)}, x_{2}]_c$ is primitive in $\toba$. Notiche that  $\bq(\alpha_1+2\alpha_2+\alpha_3,\alpha_1+2\alpha_2+\alpha_3)=1$ and $\wbq(\alpha_1+2\alpha_2+\alpha_3,\alpha_4)\neq1$, so \cite[Lemma 2.8]{AAH-triang} applies again.

\smallbreak
$\heartsuit$ Type $C_{\theta}$. Here $\wtoba_{\bq}$ is presented by the quantum Serre relations and
\begin{align} \label{eq:Ctheta-distinguished-notqsr}
[x_{(\theta-2 \, \theta)}, x_{\theta -1 \, \theta}]_c, \text { if } N=4; && \big[ [x_{(\theta-2 \, \theta)}, x_{\theta -1}]_c, x_{\theta -1} \big], \text { if } N=3.
\end{align}
As before, Proposition \ref{prop:xij=0-preNichols} gives $x_{i \theta}=0$ for any $i<\theta-1$; all the relations involving the indices $\{\theta-2, \theta -1, \theta\}$ hold in $\toba$  by Lemma \ref{lem:eminent-B3-C3} \ref{item:eminent-C3}. It remains to verify the relations involving $\{1, \dots, \theta-1\}$. Here we only have the Serre relations. But this indices span a subdiagram of type $A_{\theta-1}$, $\theta -1 \geq 3$, at a root of unity of order $3$ or $4$, so they hold by Lemma \ref{lem:eminent-A3-N>2}.
\epf

\subsection{Types $E_{6}, E_{7}$ and $E_{8}$ with $N=2$} By \cite[p. 407] {AA-diag} the distinguished pre-Nichols algebra is presented by the quantum Serre relations and
\begin{align*}
[x_{ijk}, x_j]_c=0 \, \text { if } i, j, k \text { are all different and }\widetilde{q_{ij}}, \widetilde{q_{jk}} \neq 1.
\end{align*}
\begin{lemma}\label{lem:eminent-Etheta}
If $\ba$ is of type $E_{6}, E_{7}$ or $E_{8}$ with $N=2$, then $\wtoba_{\bq}$ is eminent.
\end{lemma}
\pf
Let $\toba \in \pref (\bq)$. First we deal with the quantum Serre relations, which are always primitive. Fix $i\neq j \in \I_{\theta}$. Consider two possibilities.

\smallbreak 
$\heartsuit$ $\widetilde{q_{ij}} = 1$. In this case choose $k \in \I_{\theta}$ different from $i$ and $j$ such that $\widetilde{q_{i}} = 1$ but $\widetilde{q_{ik}} \neq 1$. We get $\bq(\alpha_{i}+\alpha_{j}, \alpha_{i}+\alpha_{j})= 1$ and $\wbq(\alpha_{i}+\alpha_{j}, \alpha_{k})\neq 1$. By \cite[Lemma 2.8]{AAH-triang}, this warranties $x_{ij}=0$ in $\toba$.

\smallbreak 
$\heartsuit$ $\widetilde{q_{ij}} \neq 1$. In this case $i$ and $j$ are consecutive vertices in a subdiagram of type $A_4$ with $N=2$. By Lemma \ref{lem:eminent-Atheta-N=2} \ref{item:Atheta-N=2-x{iij}=0} below, it follows that $x_{iij}=0$ except in the following cases: $(i,j) \in \{ (2,1), (\theta-3, \theta), (\theta-2, \theta -1)\}$. Fix such $(i,j)$, assume $x_{iij}\neq 0$ and consider $\ku x_1 \oplus \dots \oplus \ku x_{\theta} \oplus \ku x_{iij} \subset \mP(\toba)$. Then the Dynkin diagram of this braided vector space is of indefinite Cartan type. We illustrate the case $(i,j)=(\theta-3, \theta)$, the other cases being similar.
\begin{align*}
\xymatrix @C50pt 
{  
&  
&   
\overset{-1}{\underset{\theta }{\circ}} &
\overset{-1}{\underset{\theta-3 \, \theta-3 \, \theta }{\circ}} &
\\
\overset{-1}{\underset{1 }{\circ}}\ar  @{-}[r]^{-1}  &  \overset{-1}{\underset{2 }{\circ}}  \ar@{.}[r]  & 
\overset{-1}{\underset{\theta-3 }{\circ}} \ar @<0.3ex> @{-}[u]_{-1}^{} \ar
@{-}[r]^{-1} \ar @{-}[ur]_{-1}&  
\overset{-1}{\underset{\theta-2 }{\circ}}  \ar  @{-}[r]^{-1}&
\overset{-1}{\underset{\theta-1 }{\circ}}.}
\end{align*}
Thus Conjecture \ref{conj:AAH} and 
Lemma \ref{lem:subspace-of-primitives} imply $\GK \toba = \infty$.

\medbreak
Finally, fix $i,j,k$ different such that $\widetilde{q_{ij}}, \widetilde{q_{jk}} \neq 1$.  These are consecutive vertices in a suitable chosen subdiagram of type $A_4$. The Serre relations hold in $\toba$, so by Lemma \ref{lem:eminent-Atheta-N=2} \ref{item:Atheta-N=2-x2^4x1=0=x3^4x4} below we get that also $[x_{ijk}, x_j]_c=0$ in $\toba$. 
\epf

\section{On the open cases}\label{sec:open}
This section contain partial results towards  those braidings of finite Cartan type which are still open.
The detailed proofs can be found in \cite{sanmarco}.

\subsection{Type $A_2$ with $N=2$} \label{subsec:open-Cartan-A2}
\begin{lemma}\label{lem:preNichols-A2-N=2-eminent}
Assume $\ba$ is of Cartan type $A_2$ with $N=2$. Let $\toba$ be a finite GK-dimensional pre-Nichols algebra of $\bq$. The following hold:

\begin{enumerate}[leftmargin=*,label=\rm{(\alph*)}]
\item \label{item:A2-N=2-x112x221=0} if $\toba \in \pref^{\Z^2}$, then either $x_{112}=0$ or $x_{221}=0$ in $\toba$;
\item \label{item:A2-N=2-(ad4xi)xj=0} for different $i, j \in \I_2$, $(\ad_c x_i)^4x_j = 0$ in $\toba$.  \qed
\end{enumerate}
\end{lemma}

\begin{question} \label{question:A-2-N=2}
Let $\htoba_1 = \ku \langle  x_1, x_2 | x_{221}, x_{11112} \rangle$.
By Lemma \ref{lem:preNichols-A2-N=2-eminent} any $\toba \in  \pref^{\Z^2}(V)$ 
is a quotient of either $\htoba_1$ or $\htoba_2 : = \ku \langle  x_1, x_2 | x_{112}, x_{22221} \rangle$. 
Clearly $\htoba_1\simeq \htoba_2$ as algebras.  Is $\GK \htoba_1 < \infty$?
\end{question}

\subsection{Type $A_3$ with $N=2$} \label{subsec:open-typeA3-eminent}

\begin{lemma}\label{lem:eminent-A3-N=2}
Assume $\ba$ is of Cartan type $A_3$ with $N=2$. Let $\toba$ be a finite GK-dimensional pre-Nichols algebra of $\bq$. Then the following hold in $\toba$:
\begin{enumerate}[leftmargin=*,label=\rm{(\alph*)}]
\item \label{item:A3-N=2-x112=x332=0} $x_{112}=0=x_{332}$, $x_{213}\overset{\star}{=}0$, 
\item \label{item:A3-N=2-x2^41=0=x2^43} $x_{22221}=0= x_{22223}$, $x_{11113}=0=x_{33331}$,
\item \label{item:A3-N=2-x221x223=0} if $\toba \in \pref^{\Z^3}$, then at most one of $x_{113}, x_{331}, x_{221}, x_{223}$ is non-zero. \qed
\end{enumerate}
\end{lemma}

\begin{remark}
The relation $\star$ is relevant because in the tensor algebra
\begin{align}\label{eq:A3-N=2-comult-[x(13),x2]}
\begin{aligned}
\Delta ([x_{(13)}, x_2]_c)=[x_{(13)}, x_2]\ot 1 + 1 \ot [x_{(13)}, x_2] -2 q_{32} x_1 \ot x_{223} \\- 2 q_{12}^2q_{32} x_{221} \ot x_3 -2q_{12}^2q_{32} x_2 \ot x_{213} + 4 q_{12}^2q_{32} x_2^2 \ot x_{13}.
\end{aligned}
\end{align}
\end{remark}

\begin{question}\label{question:A3-N2}
By Lemma \ref{lem:eminent-A3-N=2} every $\toba \in \pref^{\Z^3}$ is 
covered by one of 
\begin{align*}
\htoba&= \ku\langle x_1, x_2, x_3 | x_{112},x_{332},x_{22221}, x_{22223},x_{11113},x_{33331},x_{213} \rangle,
\\
&\begin{aligned}
\htoba_1 &= \htoba/\langle x_{113}\rangle, &
\htoba_2 &= \htoba/\langle x_{331}\rangle, &
\htoba_3 &= \htoba/\langle x_{221}\rangle, &
\htoba_4 &= \htoba/\langle x_{223}\rangle. 
\end{aligned}
\end{align*}  
Are $\GK \htoba_1$ or $\GK \htoba_3 < \infty$? ($\htoba_1 \simeq \htoba_2$ and 
$\htoba_3 \simeq \htoba_4$ as algebras).
\end{question}

\subsection{Type $A_{\theta}$, $\theta\geq 4$ with $N=2$}
In this setting $\wtoba_{\bq}$ is presented by
\begin{align*}
x_{ij}&=0, & |i-j|&>1; &  x_{iij} &= 0, & |i-j| &= 1; &
[x_{(i i+2)}, x_{i+1}]_c &= 0, & i \in \I_{\theta-2}.
\end{align*}

\begin{lemma}\label{lem:eminent-Atheta-N=2}
Assume $\ba$ is of Cartan type $A_{\theta}$, $\theta \geq 4$, with $N=2$. The following hold in any finite GK-dimensional pre-Nichols algebra $\toba$ of $\bq$:
\begin{enumerate}[leftmargin=*,label=\rm{(\alph*)}]
\item \label{item:Atheta-N=2-x{ij}=0} $x_{ij}=0$ for any $|i-j|> 1$;
\item \label{item:Atheta-N=2-x{iij}=0} $x_{iij}=0$ for $|i-j|=1$ and $(i,j) \neq (2,1), \, (\theta-1, \theta)$;
\item \label{item:Atheta-N=2-x2^4x1=0=x3^4x4} $x_{iiiij}=0$ for  $(i,j) = (2,1), \, (\theta-1, \theta)$;
\item \label{item:Atheta-N=2-x221x334=0} if $\toba \in \pref^{\Z^{\theta}}$, then either $x_{221}=0$ or $x_{\theta-1  \theta-1  \theta}=0$;
\item \label{item:Atheta-N=2-[x(i,i+2),xi+1]=0} if $(i,j) \in \{(2,1), \, (\theta-1, \theta)\}$ and $x_{iij}=0$, then $[x_{(i-1 \ i+1)}, x_{i}]_c=0$.  \qed
\end{enumerate}
\end{lemma}

\begin{question} \label{question:A-4-N=2}
Let $\htoba_1$ denote the quotient of $T(V)$ by the relations 
\begin{align*}
x_{ij}&=0, \quad |i-j|>1; && (\ad_cx_{\theta-1})^4x_{\theta} = 0; \\
x_{iij} &= 0, \quad |i-j|=1, (i,j) \neq (\theta-1,\theta); && [x_{(13)}, x_{2}]_c=0.
\end{align*}
Similarly, define $\htoba_2$ by the relations
\begin{align*}
x_{ij}&=0, \quad |i-j|>1; && (\ad_cx_{2})^4x_1= 0; \\
x_{iij} &= 0, \quad |i-j|=1, (i,j) \neq (2,1); && [x_{(\theta-2 \ \theta)}, x_{\theta-1}]_c=0.
\end{align*}
(Clearly $\htoba_2 \simeq \htoba_1$ as algebras).  Is $\GK \htoba_1 < \infty$?
\end{question}

\subsection{Type $D_{\theta}$ with $N=2$} Here (cf. \cite[p. 404]{AA-diag}) the distinguished pre-Nichols algebra $\wtoba_{\bq}$ is presented by the quantum Serre relations 
and a  bunch of $q$-brackets coming from the several subdiagrams of type $A_3$, namely:
\begin{align}\label{eq:Dtheta-N=2-distinguished-notqsr}
[x_{(i\, i+2)}, x_{i+1}]_c, \, i\leq \theta-3; && [x_{\theta -3 \, \theta -2 \, \theta}, x_{\theta-2}]_c; && [x_{\theta \, \theta -2 \, \theta-1}, x_{\theta -2 }]_c.
\end{align}

\begin{lemma}\label{lem:eminent-D4}
Assume $\ba$ is of Cartan type $D_4$ with $N=2$. The following relations hold in any  $\toba\in \pref$:
\begin{enumerate}[leftmargin=*,label=\rm{(\alph*)}]
\item \label{item:D4-N=2-xiij=0} if $i\neq j$ and $\wbq_{ij}=-1$, then $x_{iij}=0$;
\item \label{item:D4-N=2-xij-central} if $i\neq j$ and $\wbq_{ij}=1$, then $x_{kij}=0$ for all $k \in \I_4$;
\item \label{item:D4-N=2-A3type-relations-central} if $r$ is one of the elements in \eqref{eq:Dtheta-N=2-distinguished-notqsr}, 
then $(\ad_c x_k) r=0$ for all $k \in \I_4$. \qed
\end{enumerate}
\end{lemma}

\begin{question}\label{question:D4-N2}
Let $\htoba$ denote the quotient of $T(V)$ by the relations \ref{item:D4-N=2-xiij=0}, \ref{item:D4-N=2-xij-central} and \ref{item:D4-N=2-A3type-relations-central} . Is $\GK \htoba < \infty$?
\end{question}

\begin{lemma}\label{lem:eminent-Dtheta-theta>4}
Assume $\ba$ is of Cartan type $D_{\theta}$ with $\theta >4$ and $N=2$. The following  relations hold in any  $\toba\in \pref(V)$:
\begin{enumerate}[leftmargin=*,label=\rm{(\alph*)}]
\item \label{item:Dtheta-N=2-almost} all the defining relations of $\wtoba_{\bq}$ except $x_{\theta \, \theta -1}$ and $[x_{\theta \, \theta -2 \, \theta-1}, x_{\theta -2 }]_c$;
\item \label{item:Dtheta-N=2-almost-central} the relations $x_{k \, \theta\, \theta -1}$ and $(\ad_c x_k)[x_{\theta \, \theta -2 \, \theta-1}, x_{\theta -2 }]_c$ for all $k \in \I_{\theta}$. \qed
\end{enumerate}
\end{lemma}

\begin{question}\label{question:D>4-N2}
Let $\htoba$ denote the quotient of $T(V)$ by the relations 
\begin{align*}
x_{ij}&=0, \  \wbq_{ij}=1, \  (i,j)\neq (\theta, \theta -1); &&  [x_{\theta -3 \, \theta -2 \, \theta}, x_{\theta-2}]_c=0;
\\
x_{iij}& = 0, \  \wbq_{ij}=-1;&&  [x_{(i\, i+2)}, x_{i+1}]_c=0, \ i\leq \theta-3;
\\
x_{k \, \theta\, \theta -1} &= 0, \  k \in \I_{\theta};    &&
[x_k, [x_{\theta \, \theta -2 \, \theta-1}, x_{\theta -2 }]]=0, \  k \in \I_{\theta}.
\end{align*}
Is $\GK \htoba <\infty$? 
We conjecture that $\GK \htoba= \GK\wtoba_{\bq} +2$. 
This will be treated in a subsequent paper.
\end{question}

\subsection*{Acknowledgements}
We thank Iv\'an Angiono and James Zhang for useful comments.

\end{document}